\documentclass[11 pt]{amsart}
\usepackage{graphicx, mathrsfs}
\usepackage{amssymb,amsmath,amsthm}
\usepackage{mathrsfs}
\theoremstyle{plain}
\newtheorem{theorem}{Theorem}[section]

\newtheorem{proposition}{Proposition}[section]
\newtheorem{lemma}[proposition]{Lemma}
\newtheorem{defi}{Definition}[section]

\newtheorem*{coro}{Corollary}
\newtheorem*{remark}{Remark}
\numberwithin{equation}{section} 

\newcommand{\RR}{\mathbb{R}}

\newcommand{\p}{\partial}

\newcommand{\ep}{\epsilon}

\newcommand{\CA}{\mathscr{A}}

\newcommand{\CD}{\mathcal{D}}
\newcommand{\CH}{\mathcal{H}}
\newcommand{\CM}{\mathcal{M}}
\newcommand{\TD}{\SD}

\newcommand{\D}{{\bf D}}
\newcommand{\R}{\mathcal{R}}
\newcommand{\CN}{\mathcal{N}}
\newcommand{\SD}{\mathscr{D}}
\newcommand{\SR}{\mathscr{R}}
\newcommand{\CE}{\mathcal{E}}

\newcommand{\lf}{\left}
\newcommand{\rt}{\right}

\begin{document}

\title[Water waves]{Geometry and a priori estimates for
free boundary problems of the Euler's equation}
\author[Shatah]{Jalal Shatah$^\dagger$}
\thanks{$^\dagger$ The first author is funded in part by NSF DMS
0203485.}
\address{$^\dagger$Courant Institute of Mathematical Sciences\\
251 Mercer Street\\ New York, NY 10012}
\email{shatah@cims.nyu.edu}
\author[Zeng]{Chongchun Zeng$^*$}
\address{$^*$School of Mathematics\\
Georgia Institute of Technology\\ Atlanta, GA 30332}
\email{zengch@math.gatech.edu}
\thanks{$^*$ The second author is funded in part by NSF
DMS 0627842 and the Sloan Fellowship.}
\date{}
\begin{abstract}
In this paper we derive estimates to the free boundary problem for the Euler equation with surface tension, and without surface tension provided the  Rayleigh-Taylor sign condition holds. We prove that  as the surface tension tends to zero, when the  Rayleigh-Taylor condition is satisfied, solutions converge  to the Euler flow with zero surface tension.
\end{abstract}
\maketitle

\section{Introduction}

In this paper we study free boundary problems of the Euler's
equation in vacuum:
\begin{equation} \tag{E}\begin{cases}
v_t + \nabla_v v = -\nabla p, \qquad & x\in \Omega_t
\subset \RR^n \\
\nabla \cdot v =0,& x\in \Omega_t  .
\end{cases}\end{equation}
where for every $t\in \RR$, $v(t,\cdot)$ is the velocity field of an
incompressible inviscid fluid in a  moving  domain (bounded and
connected) $\Omega_t \subset \RR^n$, $n \ge 2$,  and  $p(t,\cdot )$
is the pressure.  The boundary of the domain $\Omega_t$ moves with
the fluid velocity and the pressure at the boundary is given by the
surface tension, that is

\begin{equation}\tag{BC}\begin{cases}
\D_t=\p_t + v\cdot\nabla \quad \text{is tangent to} \quad
\bigcup\limits_{t}\Omega_t\subset \mathbb{R}^{n+1}\\
p(t,x)= \ep^2 \kappa(t,x), \qquad x\in \p \Omega_t, \quad 0\le \ep
\le 1
\end{cases}
\end{equation}
where $\kappa(t, x)$ is the mean curvature of the boundary $\p
\Omega_t$ at $x \in \p \Omega_t$, and $\D_t$ is the material
derivative. This is equivalent to saying the velocity of $\p\Omega_t$ is given by $v\cdot N$ where $N$ is the unit normal to $\p\Omega_t$ .   The case $\ep = 0$ corresponds to the zero surface
tension problem.

In the presence of surface tension we will derive  energy estimates
that bound the Sobolev norms of the velocity and the boundary.  In
addition we will show that if the Rayleigh-Taylor sign condition is
verified, then some of these bounds  are independent of
$\ep$.  In this case we conclude that as $\ep\to 0$ solutions of the
problem with  surface tension converge to solutions of the zero
surface tension problem. We do not include the effects of gravity in
(E) as it will only  contribute  lower order terms to our estimates.

The free boundary value problem for (E) has been studied intensively
by many authors. In the absence of surface tension the earliest
mathematical results on the well posedness of the water waves
problem were given by  V.~I.~Nalimov  \cite{na74} where he
considered the irrotational problem in 2~dimensions with small data
in some Sobolev space (see also H.~Yoshihara \cite{yo82}).  The
first break through in solving the well posedness for the
irrotational problem, no surface tension, for general data came in
the work of S.~J.~Wu \cite{wu97,wu99} who solved the problem in all
dimensions. For the general problem with no surface tension D. Christodoulou and H. Lindblad \cite{cl00}
were the first to obtain energy estimates   based on the  geometry of the moving domain,
assuming the Rayleigh-Taylor sign condition for  rotational flows.
H. Lindblad \cite{li05} proved existence of solutions  for the general problem.
In the absence of this condition D.~Ebin \cite{eb87} proved that the
problem is ill-posed. There is also the work of K.~ Beyer and
M.~G\"unther \cite{bg98,bg00} on well posedness which we will
comment more on at the end of the introduction.

For the problem with surface tension H.~Yoshihara \cite{yo83},
T.~Iguchi \cite{ig01} and D.~Ambrose \cite{am03} solved the
well~posedness irrotational  problem in 2~dimensions under varying
assumptions on the initial data. B.~Schweizer \cite{sc05} proved
existence for the general 3~dimensional problem.  Recently
D.~Ambrose and N.~Masmoudi  \cite{am05a} proved that as $\ep \to 0$
solutions of the 2~dimensional irrotational problem converges to
solutions of the zero surface problem by writing the equation in
terms of the arc length of the fluid boundary.

There are many other works on this problem we mention the work of
D.~Lannes \cite{la05}, T. J.~Beal, T.~Hou, and J.~Lowengrub
\cite{bhl93}, W.~Craig \cite{ca85}, G.~Shnider and E. C.~Wayne
\cite{sw02}, M.~Ogawa and A.~Tani \cite{ot02}.

During the writing of this manuscript we were informed and received
reprints of several related work by  D.~Coutand and S.~Shkoller
\cite{sc05}, D.~Ambrose  and N.~Mamoudi  \cite{am06},  and P.~Zhang
and Z.~Zhang \cite{zz06}.  In D.~Coutand and S.~Shkoller work they
proved local well posedness of the general problem using Lagrangian
coordinates.

Our results differs from those mentioned above in that we can obtain
the $\ep\to 0$ limit for the general problem.  Our approach to the
problem  is based on the well known fact that the free boundary
problem (E, BC) has a Lagrangian formulation given by
\[
I(u) = \int \int_{\Omega_0} \frac {|u_t|^2}2 dy dt -\ep^2 \int S(u)
dt,
\]
where $u(t, \cdot) \in \Gamma= \{ \Phi :\Omega_0 \to \RR^n ,
\text{volume preserving homeomorphisms}\}$, and $S(u)$ is the
surface area of $u(\p \Omega_t)$.
% When the tangent space of $\Gamma$ is endowed with the $L^2$ metric the above Lagrangian can be written as
%$$
%I(u) = \int( \frac12 {< u_t, u_t>}- \ep^2  S(u) )dt,
%$$
Critical points of $I$ satisfy
\begin{equation}\tag{E-L}
\p_tu_t + q + \ep^2S'(u) = 0,
\end{equation}
where $q$ is the Lagrange multiplier due to the constraint
$u\in\Gamma$.  Writing $v=u_t\circ u^{-1}$ and changing to Eulerian
coordinates we obtain
$$
\p_t u_t \to \D_t v, \quad q \to \nabla p_{v,v} = -\nabla
\Delta^{-1} tr(DvDv), \quad S'(u) \to J  \triangleq \nabla \kappa_\CH
$$
where $\Delta^{-1}$ is the inverse Laplacian with zero Dirichlet
data, $\kappa$ the mean curvature, and $\kappa_\CH$ is the harmonic extension of $\kappa$ into
$\Omega_t$.  Thus (E-L) is the Euler equation with the pressure $p$
given by
$$
p = p_{v,v} + \ep^2 \kappa_\CH.
$$
It is important to note that this derivation splits the pressure
into two terms, the first $p_{v,v}$ is the Lagrange multiplier, and
the second $ \kappa_\CH$ is due to surface tension. Thus these two
terms will be treated differently in the energy estimates.

Using this variational derivation, one can interpret  the Lagrange
multiplier as the second fundamental form of the manifold
$\Gamma\subset L^2(\Omega_0,\RR^n)$ and rewrite (E-L) as
\footnote[2]{Symbols in the Lagrangian description have a bar, e.g.
$\bar \SD$, while their  Eulerian counterparts do not, $\SD$.}
\begin{equation*}
\bar \SD_t u_t + \ep^2 S' (u) =0.
\end{equation*}
where $\bar \SD$ is the Riemannian connection induced on $\Gamma$ by
the embedding in $L^2$.  The above form of the equation makes it
relatively easy to identify the correct linearized problem
\begin{equation*}
\bar \SD_t^2 \bar w + \bar \SR (u_t, \bar w) u_t +\ep^2 \bar \SD^2
S(u) (\bar w) =0, \qquad \bar w(t, \cdot) \in T_{u(t, \cdot)}
\Gamma,
\end{equation*}
where $\bar \SR$ is the curvature tensor of the infinite dimensional
manifold $\Gamma\subset L^2$.  Keeping the highest order terms in
the above equation we obtain
\begin{equation} \tag{LN}
\bar \SD_t^2 \bar  w + \bar \SR_0 (v) \bar  w +\ep^2 \bar \CA \bar
w = \text{ lower order terms},
\end{equation}
where $\bar \SR_0(v)$ is a first order differential operator and
$\bar\CA$ is a third order differential operator.  In Eulerian
coordinates these terms are  given by
\begin{align*}
\SR_0(v)( w,  w) &= \int_{\p u(\Omega_0)}\ - \nabla_N
p_{v,v} |w\cdot N|^2 \; dS,\\
 \CA(u)(w,  w) &= \int_{\p u(\Omega_0)}\ |\nabla^\top w\cdot N|^2 \; dS
\end{align*}
where $N$ is the unit normal  and $\nabla^\top$ is the tangential
gradient on the boundary of $u(\Omega_0)$. Here once again we are
led in a natural way to distinguish the two problems in the following manner.
\\ 1)  For $\ep>0$  two time derivatives are  associated with $\bar \CA$ ,which is a positive semi-definite operator similar to three spatial differentiation, thus roughly speaking, $\p_t \sim (\p_x)^{\frac 32}$. Therefore one may be led to believe that the regularity of the Lagrangian coordinates given by $\p_t u=v$ is $\frac 32$ order better than $v$, which
reflects the regularizing effect of the surface tension.  However this is not true for the Lagrangian coordinates since $\CA$ is degenerate, and the regularity improvement of the $\p\Omega_t$ is geometric and is not reflected in the Lagrangian coordinates system . See Section \ref{S:counterE} for examples. Thus Eularian coordinates are more suitable to use than Lagrangian coordinates for our estimates.
\\2) For $\ep=0$
the leading term involves $\SR_0(v)$ and thus the Rayleigh-Taylor
instability may occur unless we impose the condition
\begin{equation} \tag{RT}
-\nabla_N p_{v,v}(t,x) >a >0 \quad x \in \p \Omega_t.
\end{equation}
In this case two time derivatives are  associated with $\bar\SR$ which is a positive semi-definite operator similar to
one spatial differentiation. Thus, $\p_t\sim (\p_x)^{\frac 12}$ and comments  similar to above hold on the regularity of
$\p\Omega_t$.
\\3) For $\ep>0$ one can directly obtain nonlinear estimates that depend on $\ep$ by multiplying (E-L) by $(\SD^2 S)^k S'$.
\\ 4)The control that any power  $\SR_0(v)$, with (RT) condition,
can give over vector fields is limited by the smoothness of the
boundary $\p\Omega_t$.  This fact makes the velocity  field $v$
inappropriate vector field to estimate because it is smoother than
what these operators allow.
\\5) Since  $\CA$ and  $\SR_0(v)$ are
degenerate for fields which are tangential to the boundary
$\p\Omega_t$ one needs to add the vorticity $\omega$ which controls
the rotational part of the velocity which is tangential to the
boundary.
\\These facts imply that a natural energy to control is
\[
\mathcal{E} = \int_\Omega \frac 12 |\CA^{k -1} \SD_t J|^2 + \frac
{\ep^2}2 |\CA^{k -\frac12} J|^2\; dx + \frac 12\SR_0(v)(\CA^{k -1}
J, \CA^{k -1} J) + |\omega|_{H^{3k-1} (\Omega)}^2
\]
where $J= \nabla\kappa_\CH$ is less smooth than $v$ and satisfies
\[
\bar \SD_t^2 \bar J + \bar \SR_0 (v) \bar J +\ep^2 \bar \CA \bar J =
\text{ lower order terms},
\]

In addition to the geometry of $\Gamma$ the geometry of $\p\Omega_t$
plays a crucial role in the estimates. The appearance of $\kappa$ in
the surface tension and $\nabla\kappa_\CH$ in  the energy make the
study of the geometry of the boundary as well as the study of the
harmonic extension and Dirichlet-Neumann operators on $\p\Omega_t$
central to the estimate. In using Lagrangian coordinates these
operators may be hidden but can not be avoided.

Based on these estimates one can construct an existence proof using
the following iteration method.  Since the acceleration of the
boundary is given by the surface tension plus lower order terms, in
the first step of the iteration we evolve the boundary using this
evolution. In the second step of the iteration we establish the
evolution of the velocity in the interior. This will appear in a
forthcoming paper.

Finally after this work was completed  A.~Mielke   pointed out  to
the second author that a similar geometric approach had been used by
K.~ Beyer and M.~G\"unther  to study the irrotational problem by
reducing it to the boundary\cite{bg98, bg00} . Indeed they proved
local well posedness for star shaped domains with surface tension
and studied the linearized flow for any irrotational flow. In fact
they derived the principle part of the curvature of $\Gamma$ for
irrotational  problem which of course  coincides with our $\SR_0(v)$
acting on gradient vector fields.

This paper is organized as follows. In section \ref{S:geofunc} of the paper we give the intuition behind the energy estimates by computing the geometry of $\Gamma$.  In section \ref{S:geo} we present the geometric computation of the moving boundary.  In section
 \ref{S:estiE} we present our energy estimates. In sections \ref{S:counterE} and \ref{pre} we present some examples and basic analytic and geometric calculations.

\subsection*{Notations} All notations will be defined as they are introduced.
In addition a list of symbols will be given at the end of the paper
for a quick reference.  Here we'll present some standard notations
and conventions used throughout the paper.

All constants  will be denoted by $C$ which is  a generic bound
depending only on the quantities specified in the context. We follow
the  Einstein convention where we sum upon repeated indices.

For a domain $\Omega_t$ and $x\in\p\Omega_t$ we denote by $N(t,x)$ the outward unit normal, $\Pi$ the second fundamental form where  $\Pi(w) = \nabla_w N \in T_x \p\Omega_t$ for $w\in T_x \p\Omega_t$, and $\kappa$ the mean curvature given by  the trace of $\Pi$, i.e., $\kappa = \text{tr} \Pi$.
The regularity of the domains $\Omega_t$ is characterized by the
local regularity of $\p \Omega$ as graphs. In general, an
$m$-dimensional manifold $\CM \subset \RR^n$ is said to be of class
$C^k$ or $H^s$, $s > \frac n2$, if, locally in linear frames, $\CM$
can be represented by graphs of $C^k$ or $H^s$ mappings,
respectively. For  $\p\Omega$, throughout this paper we will only
use these local graph coordinates in orthonormal frames.

%%%%%%%%%%%%%%%%%%%%%%%%%%%%%%%%%%%%%%%%%%%%%%

\section{The geometry behind the Energy}
\label{S:geofunc}

In his 1966 seminal paper \cite{ar66}, V. Arnold pointed out that the Euler equation for an incompressible inviscid fluid can be viewed as the geodesic equation on the group of volume preserving diffeomorphisms. This point of view has been adopted and developed by several authors in their work on the Euler equations on fixed domains, such as D. G. Ebin and G. Marsden~\cite{em70}, A. Shnirelman~\cite{sh85}, and Y. Brenier \cite{br99} to mention a few.  It is this point of view that  we adopt  to explain the motivation for our definition of energy.

In this section, we {\it heuristically} outline our geometric point
of view on the free boundary problems of the Euler's equation and
the intuition leading to the energy estimates in the following two
sections. Though the discussion in this section are mostly in
Lagrangian coordinates, the estimates are actually done in Eulerian
coordinates in the next two sections.

%%%%
\subsection{Lagrangian formulation of the problem} \label{SS:Lag}

One of the fundamental properties of the inviscid fluid motion is the
law of energy conservation. Multiplying the Euler's equation (E) by
$v$, integrating on $\Omega_t$, and using (BC), we obtain the
conserved energy $E_0$:
\begin{equation} \label{E:Econserv1}
E_0 = E_0(\Omega_t, v(t, \cdot))= \int_{\Omega_t} \frac {|v|^2}2 dx
+ \ep^2 \int_{\p\Omega_t} dS \triangleq \int_{\Omega_t} \frac
{|v|^2}2 dx + \ep^2 S(\p \Omega_t).
\end{equation}

The main difficulty of these problems is handling the free
boundary. A traditional way to avoid this difficulty is to
consider the Lagrangian coordinates. Let $u(t, y)$, $y \in
\Omega_0$, be the Lagrangian coordinate map solving
\begin{equation} \label{E:Lag1}
\frac {dx}{dt} = v(t, x), \qquad x(0)=y,
\end{equation}
then we have $v = u_t \circ u^{-1}$ and for any vector field $w(t, x)$, $x\in \Omega_t$, it is clear that
\begin{equation} \label{E:Dtw}
\D_t w \triangleq \p_t w + \nabla_v w = (w\circ u)_t \circ u^{-1}
\end{equation}
Therefore, the Euler's equation can be rewritten as
\begin{equation} \label{E:euler1}
u_{tt} = - (\nabla p) \circ u, \qquad u(0)= id_{\Omega_t}, \quad -
\Delta p = \text{tr}((Dv)^2), \quad p|_{\Omega_t} =\kappa,
\end{equation}
where $\kappa$ is the mean curvature of $\p \Omega_t$.

Since $v(t, \cdot)$ is divergence free, then  $u(t, \cdot)$ is volume preserving.  Let
\[
\Gamma \triangleq \{ \Phi :\Omega_0 \to \RR^n \mid \Phi \text{ is a
volume preserving homeomorphism} \}.
\]
As a manifold, the tangent space of $\Gamma$ is given by divergence
free vector fields:
\[
T_\Phi \Gamma = \{ \bar w: \Omega_0 \to \RR^n \mid \nabla \cdot w =0, \text{where } w= (\bar w \circ \Phi^{-1}) \}.
\]
%\[
%T_\Phi \Gamma = \{ \bar w: \Omega_0 \to \RR^n \mid \nabla \cdot
%(\bar w \circ \Phi^{-1})=0 \}.
%\]
For the remainder of this section we follow  the following convention: for any vector field $X \, : \Phi(\Omega_0)\to\RR^n$ its description in Lagrangian coordinates is given by $\bar X = X\circ\Phi$.
With slight abuse of notation, we also let $S(\Phi) = \int_{\p
\Phi(\Omega_0)} dS$, i.e. the surface area of $\Phi(\Omega_0)$.
Thus, the energy $E_0$ takes the following form in the Lagrangian
coordinates:
\begin{equation} \label{E:Econserv2}
E_0 = E_0 (u, u_t) = \frac 12 \int_{\Omega_0} |u_t|^2 dy + \ep^2
S(u), \qquad (u, u_t) \in T \Gamma
\end{equation}
where the volume preserving property of $u$ is used. This
conservation of energy suggests: 1) $T\Gamma$ be endowed with the
$L^2$ metric; and 2) the free boundary problem of the Euler's equation has a Lagrangian  action
\[
I(u) = \int \int_{\Omega_0} \frac {|u_t|^2}2 dx dt -\ep^2 \int S(u)
dt, \qquad u(t, \cdot) \in \Gamma.
\]

 Let $\bar \SD$ denote the covariant derivative associated with the metric on $\Gamma$, then a   critical path
  $u(t, \cdot)$ of $I$  satisfies
 \begin{equation} \label{E:critical}
\bar \SD_t u_t + \ep^2 S' (u) =0.
\end{equation}
In order to verify that the Lagrangian coordinate map $u(t, \cdot)$
of a solution of (E) and (BC) is indeed a critical path of $I$, it
is convenient to calculate $\bar \SD$ and $S'$ by viewing
$\Gamma$ as a submanifold of the Hilbert space $L^2(\Omega_0,
\RR^n)$.

%Given $\Phi \in \Gamma$, from the usual Hodge decomposition on
%$\Phi(\Omega_0)$, each vector field $\tilde w: \Omega_0 \to \RR^n$
%can be decomposed into
%\begin{equation} \label{E:decomp3}
%\tilde w= \bar w -(\nabla \psi)\circ \Phi, \qquad \bar w\in T_\Phi
%\Gamma, \quad \psi = - \Delta^{-1} (\nabla \cdot (\tilde w\circ
%\Phi^{-1})).
%\end{equation}
%Since, for any $\bar w\in T_\Phi \Gamma$ and $\psi: \Phi(\Omega_0)
%\to \RR$ satisfying $\psi|_{\p (\Phi(\Omega_0))} \equiv 0$, it holds
%\[
%<\bar w, (\nabla \psi)\circ \Phi>_{L^2 (\Omega_0)} =
%\int_{\Phi(\Omega_0)} \bar w \circ \Phi^{-1} \cdot \nabla \psi \; dx
%= 0
%\]

\noindent {\bf Computing $(T_\Phi \Gamma)^\perp$}.   For any vector field $X \; :\Phi(\Omega) \to \RR^n$ we form the Hodge decomposition

$$
X = w - \nabla \psi, \quad \psi = -\Delta^{-1}\nabla\cdot X, \quad \nabla\cdot w =0,
$$
where $\Delta^{-1}$ is the inverse Laplacian on $\Omega_t$ with zero Dirichlet data. Therefore if $\Phi\in \Gamma$,  then $\bar w = w\circ \Phi\in T_\Phi\Gamma$.  This implies
that normal space of $T_\Phi \Gamma$ at $\Phi$ is
\[
(T_\Phi \Gamma)^\perp = \{-(\nabla \psi) \circ \Phi \mid \psi|_{\p
(\Phi(\Omega_0))} \equiv 0 \}.
\]
since the Hodge decomposition is orthogonal in $L^2$ and $\Phi$ is volume preserving.

\noindent {\bf Computing $\SD_t$}. Given a path $u(t, \cdot) \in \Gamma$ and $\bar v= u_t$. Suppose
$\bar w(t, \cdot) \in T_{u(t)} \Gamma$, then the covariant
derivative $\bar \SD_t \bar w$ and the second fundamental form
$II_{u(t)} (\bar w, \bar v)$ satisfy
\[
\bar w_t = \bar \SD_t \bar w + II_{u(t)} (\bar w, \bar v), \qquad
\bar \SD_t \bar w \in T_{u(t)} \Gamma, \quad II_{u(t)} (\bar w, \bar
v) \in (T_{u(t)} \Gamma)^\perp.
\]
Let $v = u_t \circ u^{-1} = \bar v \circ u^{-1}$ and $w=\bar w \circ
u^{-1}$ which are in the Eulerian coordinates. Then from the
Hodge decomposition  we have
\begin{equation} \label{E:barCD}
\bar \SD_t \bar w = \bar w_t - II (\bar w, \bar v), \quad II (\bar
w, \bar v) = - (\nabla p_{w,v})\circ u, \qquad p_{w,v} =
-\Delta^{-1} \text{tr} (DwDv).
\end{equation}
As we do the estimates in the Eulerian coordinates, sometimes it is
more convenient to use
\begin{equation} \label{E:barCD1}
\TD_t w = (\bar \SD_t \bar w) \circ u^{-1} = \D_t w + \nabla
p_{w,v}.
\end{equation}

\noindent {\bf Computing $S'(u)$}.  By the variation of surface area  formula and for any $\bar w \in T_u \Gamma$ we have
\begin{equation} \label{E:nablaS1}
< S'(u), \bar w>_{L^2(\Omega_0)} =\int_{\p\Omega_t} \kappa
(w)^\perp \;dS = \int_{\Omega_t} \nabla \kappa_\CH
\cdot w\; dx
\end{equation}
where  $\kappa$ is the mean curvature, and $\kappa_\CH$ is its  harmonic extension.  Since $(\nabla \kappa_\CH) \circ u \in T_u \Gamma$, we
obtain
\begin{equation}\label{E:nablaS}
S'(u) = (\nabla \kappa_\CH) \circ u \triangleq J\circ u.
\end{equation}
This vector field $J$, divergence free on $\Omega_t$, is very
important for it connects the free boundary Euler's flow with the
geometry of $\p \Omega_t$ and even of $\Gamma$ as we will see later in section~\ref{SS:J}.

Combining~\eqref{E:nablaS} and~\eqref{E:barCD} with $\bar w = u_t$,
we obtain that the equation~\eqref{E:critical} for critical paths of
$I$ becomes
\begin{equation} \label{E:euler}
u_{tt} = \D_t v \circ u = (-\nabla p_{v,v} -\ep^2 J) \circ u, \qquad
v = u_t \circ u^{-1},
\end{equation}
which is equivalent to~\eqref{E:euler1}. Therefore, the free
boundary problem (E) and (BC) is a lagrangian system on  $\Gamma$  given by~\eqref{E:critical}.
 If $\ep=0$ equation~\eqref{E:critical} becomes the
geodesic equation on $\Gamma$, which is a well-known fact.

\subsection{Linearization} \label{SS:linearization}

In order to analyze the free boundary problems of the Euler's
equation, it is natural to start with the linearization. The
Lagrangian formulation provides a convenient frame work for this
purpose. From~\eqref{E:critical}, the linearized
equation is
\begin{equation} \label{E:criticalL}
\bar \SD_t^2 \bar w + \bar \SR (u_t, \bar w) u_t +\ep^2 \bar \SD^2
S(u) (\bar w) =0, \qquad \bar w(t, \cdot) \in T_{u(t, \cdot)}
\Gamma,
\end{equation}
where $\bar \SR$ is the curvature tensor of the infinite dimensional
manifold $\Gamma$.  Below we calculate $\bar \SD^2 S(u)$, which is viewed
as a linear operator on $T_u \Gamma$, and $\bar \SR$ .\\

\noindent {\bf Computing   $\bar \SD^2 S(u)$}.    Let $g(s, \cdot)$
be a geodesic  on $\Gamma$, $g(0) = u$. Let $\bar w = g_s$
and $\Omega_s = g(s, \cdot) (\Omega_0)$.
From ~\eqref{E:barCD1} we have $\bar \SD_s \bar w= (\D_s w +\nabla p_{w,w})\circ g = 0$.
Differentiating~\eqref{E:nablaS1},
\[
\bar \SD^2 S(u) (\bar w, \bar w) =\frac d{ds} \int_{\p \Omega_s} \kappa w\cdot N\; dS.
\]
%%%%%%%%%%%%%%%%%
and substitute  the expressions for $D_sN$, $D_sS$, and $D_s\kappa$ from~\eqref{E:dtn},~\eqref{E:dtdS}, and~\eqref {E:dtk2} we obtain
\begin{align*}
\bar \SD^2 S(u) (\bar w, \bar w) =&\int_{\p \Omega_s} \kappa w^\perp
(\kappa w^\perp + \CD \cdot w^\top) + \kappa \D_s w \cdot N + \kappa
w \cdot \D_s N + w^\perp
\D_s \kappa \; dS\\
=& \int_{\p \Omega_s} \kappa w^\perp (\kappa w^\perp + \CD \cdot
w^\top) - \kappa \nabla_N p_{w,w} - \kappa \nabla_{w^\top} w \cdot
N\\
&+ w^\perp\lf(-\Delta_{\p \Omega_s} w^\perp - w^\perp |\Pi|^2 + (\CD
\cdot \Pi)(w^\top)\rt) dS,
\end{align*}
where $\Pi$ is the second fundamental form of $\p \Omega_s$.

Needless to say that this is a very complicated expression for $\bar
\SD^2 S(u) (\bar w, \bar w)$. We will show that $\bar \SD^2 S(u)$ is
a differential operator and will single out its leading order part.
Let us assume that $\Omega_s$ is a sufficiently smooth domain, then
from the trace theorem,
\[
|\bar \SD^2 S(u) (\bar w, \bar w) - \int_{\p\Omega_s}\ |\nabla^\top
w^\perp|^2 + \kappa \nabla_N p_{w,w} \; dS| \le C |w|_{H^1(
\Omega_s)}^2.
\]
The third term on the left side is estimated by applying the
Divergence Theorem twice,
\[
\int_{\p \Omega_s} \kappa \nabla_N p_{w,w} dS = \int_{\Omega_s}
\nabla \kappa_\CH \cdot \nabla p_{w,w} - \kappa_\CH \text{tr} (Dw)^2
dx = \int_{\p \Omega_s} \nabla_w w \cdot \nabla \kappa_\CH dx -
\int_{\p \Omega_s} \kappa \nabla_w w \cdot N dS
\]
Therefore, we obtain
\[
|\bar \SD^2 S(u) (\bar w, \bar w) - \int_{\p\Omega_s}\ |\nabla^\top
w^\perp|^2 \; dS| \le C |w|_{L^2 (\p \Omega_s)} |w|_{H^1(\p
\Omega_s)} \le C |w|_{H^1 (\Omega_s)}^2.
\]
Much as in the derivation of~\eqref{E:nablaS1}, for a general $u \in
\Gamma$ we derive an self-adjoint operator $\bar \CA(u)$ on $T_u
\Gamma$
\[
\bar \CA (u) (\bar w) = \lf(  \nabla \CH (-\Delta_{\p u(\Omega_0)}
(w|_{\p u(\Omega_0)})^\perp)\rt) \circ u
\]
which satisfies
\[
\bar \CA(u)(\bar w, \bar w) = \int_{\p u(\Omega)}\ |\nabla^\top
w^\perp|^2 \; dS
\]
for any $\bar w \in T_u \Gamma$ and $w = \bar w \circ u^{-1}$. In
the Eulerian coordinates, $\bar \CA$ takes the form
\[
\CA (u) (w) = \nabla \CH (-\Delta_{\p u(\Omega_0)} (w|_{\p
u(\Omega_0)})^\perp), \quad \forall w: u(\Omega_0) \to \RR^n \text{
satisfying } \nabla \cdot w=0.
\]
Since $\bar \SD^2 S(u)$ is self-adjoint, then
\begin{equation} \label{E:D2S}
\bar \SD^2 S(u) = \bar \CA + \text{ at most  2nd order diff.
operators}
\end{equation}

\noindent {\bf Computing $\bar \SR$}.  In the linearized equation~\eqref{E:criticalL},
we need to calculate $\bar \SR(u) (u_t, \bar w) u_t$ for a linearized
solution $\bar w(t, \cdot)$. Therefore, we may again assume that
$u(t,\cdot)$ is a sufficiently smooth critical path of the action
$I$, thus $\bar v \triangleq v\circ u$ is smooth as
well, and study the operator $\bar \SR(u) (u_t, \cdot)u_t$ on $\bar
w$.

Here we apply a well-known formula in Riemannain geometry formally.
For any $\bar v, \bar w \in T_u \Gamma$, let $v = \bar v \circ
u^{-1}$ and $w = \bar w \circ u^{-1}$
\[
\bar \SR(u) (\bar v, \bar w)\bar v \cdot \bar w = II_u (\bar v, \bar
v) \cdot II_u (\bar w, \bar w) - II_u(\bar v, \bar w)^2 =
\int_{u(\Omega_0)} \nabla p_{v,v} \nabla p_{w,w} - |\nabla
p_{v,w}|^2 \; dx.
\]
For smooth $v$ and $w \in L^2(u(\Omega_0))$, clearly $|\nabla
p_{v,w}|_{L^2 (u(\Omega_0))} \le C |w|_{L^2(u(\Omega_t))}$. As for
the term,
\begin{align*}
\int_{u(\Omega_0)} \nabla p_{v,v} \nabla p_{w,w} \; dx =&
\int_{u(\Omega_0)} p_{v,v} \text{tr}(Dw)^2 \; dx =
\int_{u(\Omega_0)} - \nabla_w w \cdot \nabla p_{v,v}; dx \\
=&  \int_{\p u(\Omega_0)} (-\nabla_N p_{v,v}) (w^\perp)^2 \; dS +
\int_{u(\Omega_0)} D^2 p_{v,v} (w,w)\; dx.
\end{align*}
Much as in the derivation of~\eqref{E:nablaS1}, we derive an
self-adjoint operator $\bar \SR_0(v)$ on $T_u \Gamma$, depending on
$u$ and $v$,
\[
\bar \SR_0 (v) (\bar w) = \lf (\nabla \CH (- \nabla_N p_{v,v}
(w|_{\p u(\Omega_0)})^\perp)\rt) \circ u
\]
which satisfies
\[
\bar \SR_0(v)(\bar w, \bar w) = \int_{\p u(\Omega)}\ - \nabla_N
p_{v,v} w^\perp|^2 \; dS.
\]
In the Eulerian coordinates, $\bar \SR_0(v)$ takes the form
\[
\SR_0 (v) (w) = \nabla \CH (-\nabla_N p_{v,v} (w|_{\p
u(\Omega_0)})^\perp), \quad \forall w: u(\Omega_0) \to \RR^n \text{
satisfying } \nabla \cdot w=0.
\]
Therefore, in a very rough sense,
\begin{equation} \label{E:R0}
\bar \SR(u)(\bar v, \bar w)\bar v  = \bar \SR_0(\bar v) + \text{ bounded
operators}
\end{equation}
where we used the fact that $<\bar \SR(u)(\bar v, \cdot)\bar v,
\cdot>$ is self-adjoint.

\section{The  geometry of evolving  domains }
\label{S:geo}

Suppose $\Omega_t \subset \RR^n$ is a family of smooth domains with
the parameter $t$, moving with a smooth velocity vector field $v(t,
x)$, $x \in \Omega_t$. We calculate various quantities related to
the evolution of the geometry of the domain, which are essential in
the energy estimate of the free boundary problem of the Euler's
equations.

\subsection{Material derivative $\D_t$} \label{SS:mt}

For any $x_0\in \bar \Omega_{t_0}$, the particle path $x(t)$ is
the solution of the ODE:
\[
x_t = v(t, x) \qquad \qquad x(t_0) = x_0
\]
and the material derivative $\D_t = \p_t + \nabla_v$ is
differentiation along the direction of $x(t)$ in the space time
domain in $\RR^n \times \RR$. Clearly, $x(t) \in \p \Omega_t$ if
$x_0 \in \p \Omega_{t_0}$.\\

\noindent{\bf Calculations of $\D_t N$ and $\D_t S$}. At any $x_0 \in \p
\Omega_{t_0}$,  $\D_t N(t_0,x_0) \perp N(t_0,x_0)$
since $|N(t,x)|\equiv 1$. To derive $\D_t N(t_0, x_0)$, let $\tau(t) \in T_{x(t)} \p \Omega_t$ be a solution to the
linearized particle path ODE:
\[
\D_t \tau = \nabla_\tau v \qquad \qquad \tau(t_0) = \tau_0 \in T_{x_0} \p \Omega_{t_0}.
\]
At $(t_0, x_0)$,
$
\D_t N \cdot \tau_0 = \D_t (N \cdot \tau) - N \cdot \D_t \tau = - (D
v)^* (N) \cdot \tau_0
$.
Therefore, we have
\begin{equation} \label{E:dtn}
\D_t N = - ((D v)^*(N))^\top.
\end{equation}
From standard calculations for hypersurfaces
\begin{equation} \label{E:dtdS}
\D_t dS = (v^\perp \kappa + \CD\cdot v^\top) dS.
\end{equation}
%
%\noindent {\bf Calculation of $\D_t dS$.} For any $x \in
%\p\Omega_t$, let $\{y^1, \cdots, y^{n-1}\}$ be a local coordinate
%systems of a neighborhood of $x$ on $\p\Omega_t$. Let $G =
%(g_{ij})_{(n-1) \times (n-1)}$ where $g_{ij} = <\frac \p {\p y^i},
%\frac \p{\p y^j}>$, then the area element
%\[
%dS = \sqrt{ \text{det}(G)} dy^1 dy^2 \cdots dy^{n-1}.
%\]
%Thus,
%\[
%\D_t dS = \frac {\D_t (\text{det}(G))}{2 \sqrt{ \text{det}(G)}} dy^1
%dy^2 \cdots dy^{n-1} = \frac 12 \text{tr} ((\D_t G) G^{-1}) dS.
%\]
%Since the ratio between $\D_t dS$ and $dS$ is coordinate
%independent, which is the determinant of the linearization of the
%flow map generated by the vector field $v$, we may simplify our
%calculation at $x$ by requiring that $\{\frac \p {\p y^1}, \cdots,
%\frac \p {\p y^{n-1}}\}$ form an orthonormal basis at $x$. Using
%$\D_t \frac \p {\p y^i} = \nabla_{\frac \p {\p y^i}} v$, we obtain,
%at $x$,
%\[
%\text{tr} \D_t(G)= \D_t g_{ii} = 2<\nabla_{\frac \p {\p y^i}} v,
%\frac \p {\p y^i}> = 2 <\nabla_{\frac \p {\p y^i}} v^\top, \frac \p
%{\p y^i}> + 2<\nabla_{\frac \p {\p y^i}} (v^\perp N),\frac \p {\p
%y^i}>.
%\]
%It implies
%\begin{equation} \label{E:dtdS}
%\D_t dS = (v^\perp \kappa + \CD\cdot v^\top) dS.
%\end{equation}

\noindent {\bf  Covariant differentiation $\D_t^\top$.} For the
family of hypersurfaces $\p \Omega_t$ with the velocity field $v$ we define parallel transport along
the material line $x(t)$ as follows. Given a tangent vector $\tau_0 \in T_{x_0}
\p \Omega_{t_0}$, let $\tau(t)$ be the
solution of the following ODE:
\begin{equation}\label{E:ptrans}
\D_t \tau \perp T_{x(t)} \p \Omega_t \Leftrightarrow \D_t \tau =
(\nabla_\tau v \cdot N) N, \qquad \tau(t_0) = \tau_0.
\end{equation}
It is easy to verify that $\tau(t) \in T_{x(t)} \p \Omega_t$ and that  this transport preserves the
inner product.

A natural connection between $T\p \Omega_t \subset \RR^n$ for
different $t$ along the materials lines is provided by  the above parallel transport which induces the covariant
differentiation $\D_t^\top$, the projection of $\D_t$  in $\RR^n$ acting on $w\in T_{x(t)} \p \Omega_t \subset \RR^n$. This covariant differentiation induces the covariant differentiations of linear (multilinear) operators on
tensor products of $T\p\Omega_t$ and $T^*\p \Omega$, which will also denoted by $\D_t^\top$.

\noindent {\bf Calculation of $\D_t^\top \Pi$ and $\D_t \kappa$.}
Given $\tau \in T_{x_0}\p \Omega_{t_0}$, let $\tau(t)$ be its
parallel transport along the material line $x(t)$ which enable us to
compute
\[
(\D_t^\top \Pi) (\tau) = \D_t^\top ( \Pi(\tau))  = \lf( \D_t
\nabla_{\tau} N \rt)^\top = (\nabla_{\tau} \D_t N + \nabla_{[\D_t,
\tau]} N)^\top = \CD_\tau \D_t N + \Pi (([\D_t, \tau])^\top)
\]
From~\eqref{E:ptrans} and~\eqref{E:dtn}, we have
\begin{gather}\label{E:possion}
[\D_t, \tau] = \D_t \tau - \nabla_\tau (\frac {\p}{\p t} + v)
=(\nabla_\tau v \cdot N) N - \nabla_\tau v = -(\nabla_\tau v)^\top,\\
\label{E:dtpi}
(\D_t^\top \Pi) (\tau) = - \CD_\tau \lf(((D v)^*(N))^\top\rt)
-\Pi((\nabla_\tau v)^\top).
\end{gather}
%Therefore, from ~\eqref{E:dtn},
%\begin{equation}\label{E:dtpi}
%(\D_t^\top \Pi) (\tau) = - \CD_\tau \lf(((D v)^*(N))^\top\rt)
%-\Pi((\nabla_\tau v)^\top).
%\end{equation}
To calculate $\D_t \kappa$ at $(t_0, x_0)$, we take an orthonormal
frame $\{\tau_1, \ldots, \tau_{n-1}\}$ of $T_{x_0}\p \Omega_{t_0}$
and parallel transport it into an orthonormal frame  along $x(t)$.
Thus $\D_t \kappa = \D_t (\Pi(\tau_i) \cdot \tau_i) = (\D_t^\top
\Pi)(\tau_i) \cdot \tau_i$ and~\eqref{E:dtpi} give slightly different but useful forms for $\D_t \kappa$
\begin{align}
\D_t \kappa =& -\CD \cdot ((D v)^*(N))^\top - \Pi (\tau_i) \cdot
\nabla_{\tau_i} v = - \Delta_{\p \Omega} v \cdot N -
2 \Pi \cdot ((D^\top|_{T\p\Omega_t}) v) \label{E:dtk1}\\
\D_t \kappa =& -\Delta_{\p \Omega} v^\perp - v^\perp |\Pi|^2 + (\CD \cdot
\Pi)(v^\top). \label{E:dtk2}
\end{align}

\noindent {\bf Calculations of commutators involving $\D_t$.} In the
following, we will calculate the commutators of $\D_t$ with
operators $\CH$, $\CN$, and $\Delta_{\p \Omega}$, to show that they
are of lower orders.\\

\noindent{ $\bullet \quad {\displaystyle [ \D_t ,\CH]f =\Delta^{-1}
(2D v \cdot D^2 f_\CH + \nabla f_\CH \cdot \Delta v)}$.}

\noindent To start, we write the basic formula for any
function $f(t, x)$, $x \in \Omega_t$,
\begin{equation} \label{E:dtnabla}
\D_t \nabla f = \nabla \D_t f - (D v)^* (\nabla f).
\end{equation}
For the tangential gradient, using $\nabla^\top f = \nabla f -
(\nabla_N f) N$, it is straight forward to obtain
\begin{equation} \label{E:dtnablatop}
\D_t^\top \nabla^\top f = \nabla^\top \D_t f - ((D v)^* (\nabla^\top
f))^\top
\end{equation}

Let $f(t, x)$, $x \in \p \Omega_t$, be a smooth function. Recall
$f_\CH = \CH(f)$ represents the harmonic extension of $f$ into
$\Omega_t$. We have
\begin{equation}\label{E:dtdelta}
\Delta \D_t f_\CH = \D_t \Delta f_\CH + 2\nabla v \cdot D^2 f_\CH +
\nabla_{\Delta v} f_\CH = 2D v \cdot D^2 f_\CH + \nabla f_\CH \cdot
\Delta v
\end{equation}
which implies
\[
\D_t f_\CH = \CH(\D_t f) + \Delta^{-1} \Delta \D_t f_\CH = \CH (\D_t
f) + \Delta^{-1} (2D v \cdot D^2 f_\CH + \nabla f_\CH \cdot \Delta
v).
\]
Therefore we can write
\begin{equation} \label{E:dtch}
\D_t \CH(f) = \CH(\D_t f) + \Delta^{-1} (2D v \cdot D^2 f_\CH +
\nabla f_\CH \cdot \Delta v).
\end{equation}

\noindent{$\bullet\quad  {\displaystyle [\D_t, \Delta^{-1}] g =
\Delta^{-1} ( 2D v \cdot D^2 \Delta^{-1} g + \Delta v \cdot \nabla
\Delta^{-1} g) }$.}

Next, we calculate $[\D_t, \Delta^{-1}]$. Let $g(t,x)$, $x \in
\Omega_t$ be a smooth function and $\phi = \Delta^{-1} g$. From the
first half of~\eqref{E:dtdelta} where $\Delta f=0$ was not used,
\[
\D_t g = \D_t \Delta \phi = \Delta \D_t \phi - 2D v \cdot D^2 \phi-
\Delta v \cdot \nabla \phi.
\]
Since $\D_t \phi|_{\p \Omega_t} =0$, we obtain
\begin{equation} \label{E:dtdelta-1}
\D_t \Delta^{-1} g = \Delta^{-1} \D_t g + \Delta^{-1} ( 2D v \cdot
D^2 \Delta^{-1} g + \Delta v \cdot \nabla \Delta^{-1} g)
\end{equation}

\noindent{ $\bullet\quad  {\displaystyle [\D_t, \CN]f =\nabla_N
\Delta^{-1} (2D v \cdot D^2 f_\CH + \nabla f_\CH \cdot \Delta v) -
\nabla f_\CH \cdot \nabla_N v - \nabla_{\nabla^\top f} v \cdot N}$.}

To calculate the commutator of $[\D_t, \CN]$, from~\eqref{E:dtn},
\eqref{E:dtnabla} and~\eqref{E:dtch}, we have
\begin{align*}
\D_t (\nabla f_\CH \cdot N)&= \nabla_N \D_t f_\CH - \nabla f_\CH
\cdot \nabla_N v + \nabla f_\CH \cdot \D_t N\\
&= \nabla_N [\CH (\D_t f) + \Delta^{-1} (2D v \cdot D^2 f_\CH +
\nabla f_\CH \cdot \Delta v)] - \nabla f_\CH \cdot \nabla_N v -
\nabla_{\nabla^\top f} v \cdot N.
\end{align*}
Thus,
\begin{equation} \label{E:dtcn}
\D_t \CN(f) = \CN (\D_t f) + \nabla_N \Delta^{-1} (2D v \cdot D^2
f_\CH + \nabla f_\CH \cdot \Delta v) - \nabla f_\CH \cdot \nabla_N v
- \nabla_{\nabla^\top f} v \cdot N.
\end{equation}

\noindent{ $\bullet\quad  {\displaystyle[\Delta_{\p \Omega_t}, \D_t] f=
2 \CD^2 f
\cdot ((D^\top|_{T\p \Omega_t}) v)  + \nabla^\top f \cdot \Delta_{\p
\Omega_t} v  - \kappa \nabla_{\nabla^\top f} v \cdot N}$.}

In order to calculate the commutator $[\Delta_{\p \Omega_t}, \D_t]$
at $x_0 \in \p \Omega_{t_0}$, take an orthonormal frame $\{\tau_1,
\ldots, \tau_{n-1}\}$ of $T_{x_0} \p \Omega_t$. We first extend this
to an orthonormal frame to $T_x \p\Omega_{t_0}$ for all $x \in \p
\Omega_{t_0}$ close to $x_0$ by parallel transporting $\{\tau_1,
\ldots, \tau_{n-1}\}$ along geodesics on $\p\Omega_{t_0}$ starting
from $x_0$. Parallel transporting them again along the material
lines $x(t)$, we obtain an orthonormal frame $\{\tau_1, \ldots,
\tau_{n-1}\}$ of $T_x \p \Omega_t$ for all $(t,x)$ near $(t_0,
x_0)$. From the standard Riemannian geometry, this orthonormal frame
satisfies the property that, at $(t_0, x_0)$, $\CD \tau_j=0$ and
$[\tau_i, \tau_j] = \CD_{\tau_i} \tau_j - \CD_{\tau_j} \tau_i =0$,
which will be used repeatedly. For any smooth function $f(t,x)$
defined on $\p \Omega_t$, at $(t_0, x_0)$,
\begin{align*}
\D_t \Delta_{\p \Omega_t} f=& \D_t (\nabla_{\tau_j} \nabla_{\tau_j}
f - \nabla_{\CD_{\tau_j} \tau_j} f) = \nabla_{\tau_j} \D_t
\nabla_{\tau_j} f + \nabla_{[\D_t, \tau_j]} \nabla_{\tau_j} f -
\nabla_{[\D_t, \CD_{\tau_j} \tau_j]} f\\
=& \nabla_{\tau_j} \nabla_{\tau_j} \D_t f + \nabla_{\tau_j}
\nabla_{[\D_t, \tau_j]} f + \nabla_{[\D_t, \tau_j]} \nabla_{\tau_j}
f- \nabla_{[\D_t, \CD_{\tau_j} \tau_j]} f\\
=& \Delta_{\p \Omega_t} \D_t f + 2 \CD^2 f(\tau_j, [\D_t, \tau_j]) +
\nabla_{\CD_{\tau_j} [\D_t, \tau_j] -[\D_t, \CD_{\tau_j} \tau_j]} f.
\end{align*}

For any vector field $\tau(t,x) \in T_x \p\Omega_t$, it is easy to
see that $[\D_t, \tau] \in T_x\p \Omega_t$ since (a) $\tau, \frac
\p{\p t} + v \in T(\cup_t \p \Omega_t) \Rightarrow [\D_t , \tau] \in
T (\cup_t \p\Omega_t)$ and (b) $[\D_t, \tau]= \D_t \tau -\nabla_\tau
v$ does not have $\frac \p{\p t}$ component. Thus, $\CD_{\tau_j}
[\D_t, \tau_j] -[\D_t, \CD_{\tau_j} \tau_j]\in T\p \Omega_t$ and we
can drop all the normal components in its calculation. Using
$\CD_{\tau_j} \tau_j = \nabla_{\tau_j} \tau_j + \kappa N$ and $\CD
\tau_j =0$ at $(t_0, x_0)$, we obtain at $(t_0, x_0)$,
\begin{align*}
&\CD_{\tau_j} [\D_t, \tau_j] -[\D_t, \CD_{\tau_j} \tau_j] = \lf(
\nabla_{\tau_j}[\D_t, \tau_j] - \D_t (\nabla_{\tau_j} \tau_j +
\kappa N) \rt)^\top\\
=&(\nabla_{\tau_j} \D_t \tau_j - \nabla_{\tau_j}\nabla_{\tau_j}v -
\D_t \nabla_{\tau_j} \tau_j - \kappa \D_t N)^\top = -(\Delta_{\p
\Omega_t} v)^\top + \kappa ((D v)^* (N)))^\top.
\end{align*}
Therefore, from~\eqref{E:possion},
\begin{equation} \label{E:dtdeltat}
\D_t \Delta_{\p \Omega_t} f=\Delta_{\p \Omega_t} \D_t f - 2 \CD^2 f
\cdot ((D^\top|_{T\p \Omega_t}) v)  - \nabla^\top f \cdot \Delta_{\p
\Omega_t} v + \kappa \nabla_{\nabla^\top f} v \cdot N.
\end{equation}

\noindent{\bf Calculation of $\D_t^2 \kappa$.} This calculation
starts with formula~\eqref{E:dtk1}.  Since $\Pi : T\p \Omega_t \to T\p \Omega_t$ then
$\Pi \cdot D^\top|_{T\p \Omega_t} v= \Pi \cdot \nabla|_{T\p\Omega_t} v$.  Let $\{\tau_1, \ldots, \tau_{n-1}\}$ be an
orthonormal  frame which is the parallel transport of an orthonormal frame $T_{x_0} \p\Omega_{t_0}$
along the material line $x(t) \in \p \Omega_t$.
From~\eqref{E:dtn}, \eqref{E:dtpi}, \eqref{E:dtnabla}, and
\eqref{E:dtk1}, we have at $(t_0, x_0)$,
\begin{equation}\begin{split}
\D_t^2 \kappa =& - \D_t \Delta_{\p \Omega_t} v \cdot N - \Delta_{\p
\Omega_t} v \cdot \D_t N - 2 (\D^\top_t (\Pi (\tau_i))) \cdot
\nabla_{\tau_i} v - 2\Pi
(\tau_i) \cdot \D_t (\nabla_{\tau_i} v)\\
= & - \D_t \Delta_{\p \Omega_t} v \cdot N + \Delta_{\p \Omega_t} v
\cdot (D v)^*(N)^\top + 2 \CD_{\tau_i} \lf(((D v)^*(N))^\top\rt)
\cdot \nabla_{\tau_i} v  \\
&+ 2\Pi((\nabla_{\tau_i} v)^\top) \cdot \nabla_{\tau_i} v -
2\Pi(\tau_i) \cdot \nabla_{\D_t \tau_i} v - 2\Pi(\tau_i) \cdot
\nabla_{\tau_i} \D_t v + 2\Pi(\tau_i) \cdot (D v)^2 (\tau_i)  \\
=& - \D_t \Delta_{\p \Omega_t} v \cdot N - 2 \Pi \cdot (D^\top|_{T
\p \Omega_{t_0}} \D_t v) + \Delta_{\p \Omega_t} v \cdot (D
v)^*(N)^\top + 2[ \CD \lf(((D v)^*(N))^\top\rt)   \\
&+ \Pi((D^\top|_{T \p \Omega_{t_0}} v)^\top)] \cdot (D^\top|_{T\p
\Omega_{t_0}} v) + 2\Pi \cdot ((Dv)^2|_{T \p \Omega_t})^\top.
\label{E:dttk1}\end{split}
\end{equation}

To compute $\D_t \Delta_{\p \Omega_t} v \cdot N $ from~\eqref{E:dtdeltat} we need the general formula
\[
\CD^2 f(\tau, \tau') = D^2 f(\tau, \tau') - (\Pi(\tau) \cdot \tau')
\nabla_N f.
\]
for any $\tau, \tau' \in T_{x_0}\p \Omega_{t_0}$. Therefore,
\begin{equation}\begin{split}
- \D_t \Delta_{\p \Omega_t} v \cdot N = & - N \cdot \Delta_{\p
\Omega_t} \D_t v + 2 N \cdot D^2 v (\tau_i, (\nabla_{\tau_i}v)^\top)
- 2 (\nabla_N v \cdot
N) (\Pi \cdot (D^\top|_{T\p \Omega_{t_0}} v)) \\
&+ N \cdot \nabla v ((\Delta_{\p \Omega_t} v)^\top) - \kappa |
(\nabla v)^*(N)^\top|^2. \label{E:dttk2}
\end{split}\end{equation}
When $v$ and $\Omega_t$ satisfy the Euler's equation, the expression
for $\D_t^2 \kappa$ can be written as
\begin{equation} \label{E:dttk3}
\D_t^2 \kappa = -N \cdot \Delta_{\p \Omega_t} \D_t v + 2 \ep^2
\Pi \cdot (D^\top|_{T \p \Omega} J) +   r
\end{equation}
where we signaled out the important terms in the above equation
\subsection{An important vector field for the water wave problem}
\label{SS:J}

Since $J=\nabla \kappa_\CH$ appears in the Euler's equation as a force generated by the
surface tension and its regularity is closely related to that of $\p
\Omega_t$, we will study the temporal evolution of $J$ for
the rest of this section.

\noindent{\bf Computing $\D_t J$ and $\D_{t}^2J$.}
From~\eqref{E:dtk1}, \eqref{E:dtnabla}, \eqref{E:dtch} and the
definition of $J$,
\begin{equation}\begin{split}
\D_t J = &\nabla \D_t \kappa_\CH - (D v)^* J = \nabla \CH (\D_t
\kappa) + \nabla \Delta^{-1} [ 2 D v \cdot D J + J \cdot \Delta v]
- (D v)^* J \\
=& \nabla \CH (\Delta_{\p \Omega_t} v \cdot N) - 2 \nabla \CH(\Pi
\cdot ((D^\top|_{T\p\Omega_t}) v))  \\
& + \nabla \Delta^{-1} [ 2 D v \cdot D J + J \cdot \Delta v]
- (D v)^* J. \label{E:dtJ1}\end{split}
\end{equation}
From~\eqref{E:dtk2}, a slightly different way to write $\D_tJ$ is
\begin{equation} \label{E:dtJ2} \begin{split}
\D_t J = -\nabla \CH(\Delta_{\p \Omega_t} v^\perp) + \nabla \CH[
-v^\perp |\Pi|^2 + (\CD \cdot \Pi)(v^\top)]\\ + \nabla \Delta^{-1} [
2 D v \cdot D J + J \cdot \Delta v] - (D v)^*J.
\end{split}\end{equation}
Generally, when the surface tension is of order $O(1)$, it is
sufficient to consider $\D_t J$. However, when there is no surface
tension or the surface tension converges to $0$, we have to
calculate $\D_t^2 J$. Differentiating~\eqref{E:dtJ1}, we obtain
\begin{equation}\begin{split}
\D_t^2 J = & \D_t (\nabla \D_t \kappa_\CH - (D v)^* J) \\
=& \nabla \D_t^2 \kappa_\CH - (Dv)^* \nabla \D_t \kappa_\CH
- (Dv)^* \D_t J - (\D_t Dv)^* J \\
=& \nabla \D_t^2 \kappa_\CH - 2(Dv)^* \D_t J - ((Dv)^2)^* J- (D
\D_t v - (Dv)^2)^* J. \label{E:dttJ1}\end{split}
\end{equation}
Since $\kappa_\CH$ is harmonic, from~\eqref{E:dtch}
and~\eqref{E:dtJ1},
\begin{equation}\begin{split}
\D_t^2 \kappa_\CH = &\D_t \lf[\CH (\D_t \kappa) + \Delta^{-1} ( 2 D
v
\cdot D J + J \cdot \Delta v) \rt] \\
= & \CH(\D_t^2 \kappa) + \Delta^{-1} ( 2 D v \cdot D^2  + \Delta v
\cdot \nabla) \CH(\D_t \kappa) + \D_t \Delta^{-1} ( 2 D v \cdot D J
+ J \cdot \Delta v) \\
=& \CH(\D_t^2 \kappa) + \Delta^{-1} ( 2 D v \cdot D + \Delta v
\cdot) [\D_t J - \nabla \Delta^{-1} ( 2 D v \cdot D J + J \cdot
\Delta v) - (D v)^* J] \\
& + \D_t \Delta^{-1} ( 2 D v \cdot D J + J \cdot \Delta v)
\label{E:dt2chk}\end{split}
\end{equation}
A more explicit expression of $\CH(\D_t^2 \kappa)$ can be derived
from~\eqref{E:dttk3}.  The $\D_t\Delta^{-1}$ is another term in the above equation that can be explicitly calculated from~\eqref{E:dtdelta-1} and~\eqref{E:dtdelta} to write
\begin{equation}\begin{split}
\D_t \Delta^{-1} &( 2 D v \cdot D J + J \cdot \Delta v)\\
=& \Delta^{-1} \D_t ( 2 D v \cdot D J + J \cdot \Delta v) +
\Delta^{-1} [ (2Dv \cdot D^2 + \Delta v \cdot \nabla)
\Delta^{-1} ( 2 D v \cdot D J + J \cdot \Delta v)]\\
=& \Delta^{-1} [ 2 (D\D_t v -  (D v)^2) \cdot DJ + 2 Dv \cdot (D
\D_t J - D J Dv) + \D_t J \cdot \Delta v \\
&+ J \cdot (\Delta \D_t v - \nabla_{\Delta v} v - 2 \p_i v^j \p_{ij}
v)  ] \\
&+ \Delta^{-1} [ (2Dv \cdot D^2 + \Delta v \cdot \nabla) \Delta^{-1}
( 2 D v \cdot D J + J \cdot \Delta v)] \label{E:dt2chk1}
\end{split}
\end{equation}

%%%%%%%%%%%%%%%%%%%%%%%%%%Insert 1

\subsection*{\bf Computing $ \SD_t J$ and $\SD_{tt} J$ with divergence free $v$.}
In the rest of this section,
assume $\nabla \cdot v=0$. Given any vector field $w$ defined on
$\Omega_t$ with $\nabla \cdot w=0$, let $\TD_t w$ denote the
divergence free part of $\D_t w$. It is easy to calculate that
\begin{equation} \label{E:TDII}
\TD_t w = \D_t w + \nabla p_{v,w}, \qquad  -\Delta p_{v,w} =
\text{tr}(D v D w), \quad p_{v,w}|_{\p \Omega_t} =0.
\end{equation}

As $J$ is divergence free, we will decompose material derivatives of
$J$ into the divergence parts and gradient parts, i.e. we consider
$\TD_t J$ and $\TD_{tt} J$, the covariant derivatives defined
in~\eqref{E:TDII}. Then we have
\begin{equation}\label{E:tdtJ}
\TD_t J= \D_t J + \nabla p_{v, J}.
\end{equation}
For the second order derivative,
\begin{equation} \label{E:tdt2J1}
\TD_t^2 J = \D_t \TD_t J + \nabla p_{v, \TD_t J} = \D_t^2 J + \D_t
\nabla p_{v, J} + \nabla p_{v, \TD_t J}.
\end{equation}
where $p_{v, J}$ and $p_{v, \TD_t J}$ are defined as
in~\eqref{E:TDII}. Using~\eqref{E:dtnabla}, \eqref{E:dtdelta-1}

\begin{equation}\begin{split}
\D_t \nabla p_{v, J} =& \nabla \D_t p_{v, J} - (Dv)^* \nabla p_{v,
J} = \nabla \D_t \Delta^{-1} (Dv \cdot D J) - (Dv)^* \nabla p_{v,
J} \\
=& \nabla \Delta^{-1} \D_t (Dv \cdot DJ) + \nabla \Delta^{-1} (2 Dv
\cdot D^2 + \Delta v \cdot D) p_{v, J} - (Dv)^* \nabla p_{v,
J} \\
=& \nabla \Delta^{-1} [(D\D_t v -  (D v)^2) \cdot DJ + Dv \cdot (D
\D_t J - D J Dv)] \\
&+ \nabla \Delta^{-1} (2 Dv \cdot D^2 + \Delta v \cdot D) p_{v, J} -
(Dv)^* \nabla p_{v, J}. \label{E:tdt2J2}\end{split}
\end{equation}
Using the above calculations, we will show in Lemma~\ref{L:J} that
$J$ satisfies the linearized Euler's equation with lower order
terms. The estimates on $J$ and $\D_t J$ from the linearized Euler's
Equation will imply the estimates on the geometry of the moving
domain and the velocity fields.

%%%%%%%%%%%%%%%%%%%%%%%%%%%%%%%%%%%%%%%%%%%%%%

\section{Main Results} \label{S:estiE}

In this section, we will derive local energy estimates and prove
convergence theorems. We show that solutions of (E) with boundary
condition (BC) are locally bounded
\begin{equation} \label{E:assump}
v(t, \cdot) \in H^{3k} (\Omega_t) \quad \text{ and } \quad \p
\Omega_t \in H^{s_0}, \qquad s_0=3k \text{ or } 3k +1 \text{ for }
\ep=0 \text{ or } >0
\end{equation}
where $k$ is an integer satisfying $3k> \frac n2 +1$ (equivalently
$3k \ge \frac n2 + \frac 32$). When $\ep>0$, this estimate is
obtained without any additional assumption and it may depend on
$\ep$. To derive a priori estimates independent of $\ep$, we assume
the Rayleigh-Taylor sign condition (RT): $-\nabla_N p_{v,v}(t,x)
>a >0 \quad x \in \p \Omega_t$ for some constant $a$. \\

\noindent{\bf Definition of the energies and statements of the
theorems}. The conserved energy of  the Euler's equation is given by
\[ E_0 = E_0 (\Omega,
v) = \int_\Omega \frac 12 |v|^2 dx + \ep^2 S(\Omega),
\]
where $S(\Omega) = \int_{\p \Omega} dS$ is the surface area. Higher
order energies are based on the linearized Euler flow and thus
involve the differential operators $\SD_t$, $\CA$, and $\SR_0(v)$.

Recall that, for any  vector field $w \in H^s (\Omega)$ with $\nabla
\cdot w=0$, $\CA$ is given by
\begin{equation}\label{E:CA0}
\CA(w) = - \nabla \CH \Delta_{\p \Omega}\, w^\perp.
\end{equation}
$\CA$ is a semi-positive definite self-adjoint third order
differential operator if $w|_{\p \Omega}$ is not smoother than $N
\in H^{s_0-1} (\p\Omega)$. In fact, it is positive definite acting
on the irrotational part $w_{ir}$ of $w$, i.e. $w_{ir}=\nabla \CH
\CN^{-1} w^\perp$, see~\eqref{E:decom2} for details.

Also recall that $\SR_0(v)$, which depends on $\Omega$ as well as on a
vector field $v \in H^{3k} (\Omega)$ with $\nabla \cdot v =0$, is
given by
\begin{equation} \label{E:TR0}
\SR_0(v) (w) = \nabla \CH ((-\nabla_N\, p_{v,v}) w^\perp) =  \nabla \CH((-\nabla_N \Delta^{-1} tr(DvDv)) w^\perp) .
\end{equation}
By Lemma~\ref{L:Delta}, $\nabla_N p_{v,v} \in H^{3k- \frac 32} (\p
\Omega)$ and therefore $\SR_0(v)$ is a first order self-adjoint
differential operator if $w$ is not smoother than $\nabla p_{v,v}$.
Under the sign assumption (RT), $\SR_0(v)$ is semi-positive
definite, and like $\CA$, it is positive definite on the
irrotational part $w_{ir}$ of $w$.

Let $\omega_v: \RR^n \to \RR^n$, often simply written as $\omega$
for short, represent the curl or vorticity of a vector field $v$
defined on $\Omega$, i.e.
\[
\omega (X) \cdot Y = \nabla_X v \cdot Y - \nabla_Y v \cdot X
\]
for any vector $X, Y \in \RR^n$. Viewing $\omega$ as a matrix, its
entries are $\omega_i^j = \omega(\frac \p{\p x^i}) \cdot \frac \p{\p
x^j} = \p_i v^j - \p_j v^i$.

\begin{defi} \label{E:energy}
For any domain $\Omega$ in $H^{s_0}$, $s_0 =3k$ if $\ep=0$ or
$s_0=3k+1$ if $\ep>0$, and any vector field $v \in H^{3k}(\Omega)$
with $\nabla \cdot v=0$, define the energies $E(\Omega, v)$ and
$E_{RT} (\Omega, v)$ , often written as $E$ and $E_{RT}$ for short,
\[
E = \int_\Omega \frac 12 |\CA^{k -1} \SD_t J|^2 + \frac {\ep^2}2
|\CA^{k -\frac12} J|^2\; dx + |\omega|_{H^{3k-1} (\Omega)}^2, \qquad
E_{RT} = \int_\Omega \frac 12 \SR_0(v) \CA^{k -1} J \cdot \CA^{k -1}
J dx
\]
where $\SD_t$ is the divergence free part of $\D_t J$ defined
in~\eqref{E:TDII},
\[
\TD_t J = \D_t J + \nabla p_{v, J} = \D_t J - \nabla \Delta^{-1}
\text{tr} (DvDJ)
\]
Set $\CE =E+ E_{RT}$.
\end{defi}
To estimate terms in the energy we need to consider  the following type of
$H^{s_0}$ neighborhoods of $\Omega_*$,  a bounded connected domain in $\RR^n$,  which are  bounded
in $H^s$ for some $s \ge s_0$.
\begin{defi}\label{D:domainnbd} Let $\Lambda= \Lambda (\Omega_*,s_0,
s, L, \delta)$ be the collection of all domains $\Omega$ satisfying
\begin{enumerate}
\item [(A1)] there exists a diffeomorphism $F: \p \Omega_* \to \p
\Omega \subset \RR^n$, so that $|F - id_{\p \Omega_*}|_{H^{s_0} (\p
\Omega_*)} < \delta$;
 \item [(A2)]  the mean curvature $\kappa$ of $\p \Omega$ satisfies
$|\kappa|_{H^{s-2}(\p \Omega)} < L$.
\end{enumerate}
\end{defi}

%Let $\Omega_* \subset \RR^n$ be a bounded connected domain
%with $\p \Omega_* \in H^{s_0}$. We conside

Fix $0 < \delta \ll 1$ and let $L_0 =1+ 2 |\kappa(0, \cdot)|_{H^{3k -
\frac 52}(\p \Omega)}$ and  $\Lambda_0 \triangleq \Lambda
(\Omega_0, 3k - \frac 12, 3k-\frac 12, L_0, \delta)$.   By
Lemma~\ref{L:Delta}, and equations  \eqref{E:trace}  and~\eqref{E:vperp}  we have
\begin{align}
&|\CA|_{L(H^s(\Omega), H^{s-3}(\Omega))} \le C, && s \in [3,
3k-1]\\
&|\SR_0(v)|_{L(H^s (\Omega), H^{s-1} (\Omega))}\le C |\nabla_N \,
p_{v,v}|_{H^{3k -\frac32} (\p \Omega)} \le C |v|_{H^{3k}(\Omega)}^2,
\qquad && s \in [1, 3k -1]\\
&|p_{v,v}|_{H^{3k}(\Omega)} \le C |(Dv)^2|_{H^{3k-2} (\Omega)} \le C
|v|_{H^{3k-\frac 58} (\Omega)}^2 \qquad \text{by Sobolev
inequalities} \label{E:pressure1}
\end{align}
where $C$ is uniform in $\Omega \in \Lambda_0$. Here used the fact
$3k \ge \frac n2 + \frac 32$.  The norm $H^{3k-\frac 58}$ in \eqref{E:pressure1} is chosen for convenience; any norm $H^{3k-1+\alpha}$ would work with $\alpha >\frac 14$.  The next proposition gives bounds on the velocity and mean curvature in terms of these energies.
\begin{proposition} \label{P:energy}
For $\Omega \in \Lambda_0$ with $\p \Omega \in H^{s_0}$, we have
\[
\ep^2 |\kappa|_{H^{3k -1}(\p \Omega)}^2 \le 3E + C_0\ep^2, \qquad
|v|_{H^{3k} (\Omega)}^2 \le C_0 (E+ E_0)
\]
and, if we also assume (RT),
\[
|\kappa|^2_{H^{3k-2}(\p \Omega)} \le C_* E_{RT} + C_0,
\]
for some constant $C_*, C_0>0$. $C_*$ depends only on $a$ in
assumption (RT) and $C_0$ depends only on the set $\Lambda_0$.
\end{proposition}

The proof of this proposition will be given below. Using this result
we will prove the following three theorems. The first theorem holds
when there is surface tension which makes the regularity of
$\p\Omega_t$ better (in $H^{3k+1}$) but the bound on
$|\kappa|_{H^{3k-1}(\p\Omega_t)}$ depends on $\ep$, the strength of
of the surface tension.

\begin{theorem} \label{T:STenergyE}
Assume $\ep>0$ and fix $\delta>0$ sufficiently small. Then there
exists $L_\ep>0$ such that, if a solution of (E) and (BC) is given
by $\Omega_t$ with $\p \Omega_t \in H^{3k+1}$ and $v(t, \cdot)\in C(H^{3k}(\Omega_t))$,
then there exists $t^*>0$,
depending only on $|v(0, \cdot)|_{H^{3k} (\Omega_t)}$, $L_\ep$, and
the set $\Lambda _0$, such that, for all $t \in [0, t^*]$,
\begin{equation}\begin{split}
&\Omega_t \in \Lambda_0 \qquad \text{ and } \qquad
|\kappa|_{H^{3k-1}(\p \Omega_t)} \le L_\ep, \\
&E(\Omega_t, v(t, \cdot)) \le 2 E(\Omega_0, v(0, \cdot)) + C_1 +
\int_0^t P_\ep (E_0, E(\Omega_{t'}, v(t', \cdot))) \, dt'
\label{E:STenergyE}
\end{split}
\end{equation}
where $P_\ep(\cdot)$ is a polynomial of positive coefficients
determined only by $\ep$ and the set $\Lambda_0$ and $C_\ep$ is an
constant determined only by $\ep$, $|v(0, \cdot)|_{H^{3k
-\frac32}(\Omega_0)}$, and the set $\Lambda_0$.
\end{theorem}

Since the domain is evolving, the above continuity assumption of $v$
in $t$ means that there exists an extension of $v$ to $[0, T] \times
\RR^n$ which is continuous in $H^{3k}(\RR^n)$. The second theorem
holds under the assumption (RT) and the estimates are uniform in
$\ep \in [0,1]$. As it does not take the advantage of the surface
tension even if it is present, the bound on the regularity of $\p
\Omega_t$ is only in $H^{3k}$.

\begin{theorem} \label{T:energyE}
Assume $\ep\in [0,1]$ and (RT) holds. Fix sufficiently small
$\delta>0$. There exists $L>0$ such that, if a solution of (E) and
(BC) is given by $\Omega_t$ with $\p \Omega_t \in H^{3k}$ and $v(t, \cdot)\in C(H^{3k}(\Omega_t))$, then there exists $t^*>0$,
depending only on $|v(0, \cdot)|_{H^{3k} (\Omega_t)}$, $L$, and the
set $\Lambda_0$, such that, for all $t \in [0, t^*]$,
\begin{equation}\begin{split}
&\Omega_t \in \Lambda_0 \qquad \text{ and } \qquad
|\kappa|_{H^{3k-2} (\p \Omega_t)} \le L, \\
&\CE(\Omega_t, v(t, \cdot)) \le \CE(\Omega_0, v(0, \cdot)) +
\int_0^t P(\CE(\Omega_{t'}, v(t', \cdot))) \, dt' \label{E:energyE}
\end{split}
\end{equation}
where $P(\cdot)$ is a polynomial of positive coefficients uniform in
$\ep$, determined by the set $\Lambda_0$.
\end{theorem}

An immediate consequence of the above theorem is convergence of solution as
the surface tension approaches $0$.

\begin{theorem} \label{T:VST}
Assume (RT) holds. Fix the initial data $\p \Omega_0 \in H^{3k +1}$
and $v(0, \cdot) \in H^{3k} (\Omega_0)$. As $\ep \to 0$, subject to
a subsequence, the solution of (E) and (BC) with vanishing surface
tension converges to a solution of (E) and (BC) for $\ep=0$ weakly
in the space of $\p \Omega_t \in H^{3k}$ and $v(t, \cdot) \in
H^{3k}$.
\end{theorem}

The above convergence of $\p\Omega_t$ is in the sense of local
coordinates and the convergence of $v$ can be obtained by using the
Lagrangian coordinates $u(t,y)$ which is also in $H^{3k}$. We also
observe that  the neighborhood $\Lambda_0$ of the domains
$\Omega$ does not have to be centered at $\Omega_0$. Thus, since the
constants involved in the energy estimates only depend on the
neighborhoods and the norm of the initial velocity, these estimates
provide a basis for a continuation argument local in time.\\

\noindent {\bf Proof of Proposition~\ref{P:energy}.} From the
definition of $\CA$ and $\SR_0(v)$, it is easy to obtain
\begin{align*}
&\ep^2 |\nabla^\top \CN \lf(-\Delta_{\p \Omega} \CN \rt)^{k-1}
\kappa|_{L^2(\p \Omega)}^2 = \ep^2 |\CA^{k-\frac 12}
J|_{L^2(\Omega)} \le 2 E, \\
&\int_{\p \Omega} -\nabla_N p_{v,v} |\CN (-\Delta_{\p \Omega} \CN)^{
k -1} \kappa|^2 \; dS = <\SR_0(v) \CA^{k -1} J,
\CA^{k-1}J>_{L^2(\Omega)}^2 = 2E_{RT}.
\end{align*}
To estimate $\kappa$ in either $H^{3k-2}(\p \Omega)$ or $H^{3k-1}
(\p \Omega)$, it is sufficient to use the estimates on functions and
operators defined on $\p \Omega$ considered as in $H^{3k-\frac 12}$
only. Therefore, the inequality for $\kappa$ in
Proposition~\ref{P:energy} follows from (RT) and the fact that $\CN$ behaves as a first order derivative (Theorem~\ref{T:DeltaN},
 \eqref{E:deltacn3}, and \eqref{E:deltacn2}).

To estimate $v$, it is easy to calculate that
\begin{equation} \label{E:deltav}
\Delta v^i = \p_j \omega_j^i
\end{equation}
which is part of the energy. Therefore, we only need to show that
some boundary data of $v$ is controlled by  $E$ and   the conserved energy $E_0$. This
boundary data of $v$ turns out to be $\nabla_N \, v$. Since
$\nabla_N \, v = (D v)^*(N) + \omega (N)$ and $\omega$ is controlled
by $E$, it suffices to
estimate $\nu = (D v)^*(N)$. \\

\noindent {\it Step 1. Tangential curl $\omega_\nu^\top$ of
$\nu^\top$.} Let $\nu^\top$ be the tangential component of $\nu$ and
$\omega_\nu^\top(x): T_x\p\Omega \to T_x\p \Omega$ be defined as
\[
\omega_\nu^\top (x) (X) \cdot Y = \CD_X \, \nu^\top \cdot Y -
\CD_Y\, \nu^\top \cdot X = \nabla_X \, \nu^\top \cdot Y - \nabla_Y
\, \nu^\top \cdot X
\]
for any $x \in \p \Omega$ and $X, Y \in T_x \p \Omega$. To obtain a
more explicit form of $\omega_\nu^\top$, let $X$ and $Y$ be extended
to tangent vector fields on a neighborhood of $x$ on $\p \Omega$ by
parallel transport along geodesics on $\p \Omega$ emitting from $x$.
From the definition of $\nu$, we have
\begin{equation}\begin{split}
\omega_\nu^\top (X) \cdot Y = & \nabla_X ( \nabla_Y \, v \cdot N) -
\CD_X Y \cdot \nu - \nabla_Y ( \nabla_X \, v \cdot N) + \CD_Y X\cdot
\nu \\
= &\Pi(X) \cdot \nabla_Y \, v - \Pi(Y) \cdot \nabla_X \, v.
\label{E:tcurl1}
\end{split}
\end{equation}
Therefore, by Sobolev inequalities, there exists $C>0$ uniform in
$\Omega \in \Lambda_0$ so that
\[
|\omega_\nu^\top|_{H^{3k - \frac 52} (\p \Omega)} \le C|\Pi \circ
D^\top v|_{H^{3k - \frac 52} (\p \Omega)} \le C |v|_{H^{3k -\frac
18} (\Omega)}
\]
since $3k \ge \frac n2 +\frac32$. Again here the norm $H^{3k-\frac18}$ is chosen to illustrate
that the term is lower order. In fact any $H^{3k -\alpha}$ with $0 < \alpha < \frac 14$ works.

\noindent {\it Step 2. Divergence $\CD \cdot \nu^\top$.} At any $x_0
\in \p \Omega$, let $\{X_1, \ldots, X_{n-1}\}$ be an orthonormal
frame of $T_{x_0} \p \Omega$. We extend them to orthonormal frames
of $T_x \p \Omega$ at each $x$ in a neighborhood of $x_0$ in $\p
\Omega$ by parallel transport along geodesics on $\p \Omega$
emitting from $x_0$. At $x_0$,
\begin{equation}\label{E:cddotv1}
\CD \cdot \nu^\top = \CD_{X_i}\, \nu^\top \cdot X_i = \CD_{X_i}
(\nabla_{X_i}\, v \cdot N) = \Delta_{\p \Omega} \, v \cdot N +
(D^\top|_{T \p\Omega_t}) v \cdot \Pi.
\end{equation}
To control the first term on the right side, we use~\eqref{E:tdtJ}
and~\eqref{E:dtJ1} to obtain
\[
|\SD_t J + \nabla \CH (\Delta_{\p \Omega} \, v \cdot N)|_{H^{3k-3}
(\Omega)} \le C |v|_{H^{3k-\frac 18 } (\Omega)}
\]
where we also used Lemma~\ref{L:Delta} for the estimates. To get boundary estimates, note that
$\SD_t J \cdot N$ and $N \cdot \nabla \CH (\Delta_{\p \Omega} \, v
\cdot N)$ are well defined in $H^{3k -\frac 72} (\p \Omega)$ even if
$k=1$ since they are divergence free. and thus combining the above inequality
with the identity for $\CD \cdot \nu^\top$ \eqref{E:cddotv1}, we obtain
\[
|\CN (\CD \cdot \nu^\top) + (\SD_t J)^\perp|_{H^{3k - \frac72} (\p
\Omega)} \le C |v|_{H^{3k -\frac 18} (\Omega)} \quad \Rightarrow
\quad |\CD \cdot \nu|_{H^{3k-\frac52} (\p \Omega)}^2 \le C (E +
|v|^2_{H^{3k-\frac 18} (\Omega)}).
\]

\noindent {\it Step 3. Control of $\nu^\top$.} Using the same frame
as in Step 2, at $x_0$,
\begin{align*}
\Delta_{\p \Omega} \nu^\top \cdot X_i = & \CD_{X_j} \CD_{X_j}\,
\nu^\top \cdot X_i = \nabla_{X_j} (\CD_{X_j}\, \nu^\top \cdot X_i) =
\nabla_{X_j} (\CD_{X_i}\, \nu^\top \cdot X_j + \omega_\nu^\top (X_j)
\cdot X_i)\\
=& \nabla_{X_i} (\CD \cdot \nu^\top) + \R(X_i, X_j) \nu^\top \cdot
X_j + (\CD_{X_j}\, \omega_\nu^\top) (X_j) \cdot X_i.
\end{align*}
One can also write
\[
\Delta_{\p \Omega}\, \nu^\top = \nabla^\top (\CD \cdot \nu^\top) +
\text{Ric}((\nabla v)^*(N)^\top) +(\CD_{X_j}\, \omega_\nu^\top)
(X_j)
\]
where Ric is the Ricci curvature of $\p \Omega$. From the estimate
on $\omega_\nu^\top$ and $\CD \cdot \nu^\top$,
\[
|\nu^\top|^2_{H^{3k-\frac 32} (\p \Omega)} \le C (E +
|v|^2_{H^{3k-\frac 18} (\Omega)})
\]
with a uniform constant $C>0$. \\

\noindent {\it Step 4. Normal component $\nu^\perp = \nabla_N\, v
\cdot N$ of $\nu$.} This will be estimated by calculating the
divergence of $\nu$ in two ways. Recall $N_\CH = (N_\CH^1, \ldots,
N_\CH^n)$ denotes the harmonic extension of $N$ into $\Omega$. Let
$\nu$ also be extended to $(D v)^*(N_\CH)$. Near any $x_0 \in \p
\Omega$, let $\{X_1, \ldots, X_{n-1}\}$ be the orthonormal frame of
$T_x \p \Omega$ constructed above and let $X_n = N_\CH$. On one
hand, at $x_0$,
\begin{align*}
\nabla \cdot \nu = &\nabla_{X_i} ((D v)^* (N_\CH)) \cdot X_i =
\nabla_{X_i} (\nabla_{X_i} \, v \cdot N_\CH) - N \cdot D v
(\nabla_{X_i} X_i)\\
=& (\nabla_{X_i} \omega)(X_i) \cdot N + D v \cdot D N_\CH
\end{align*}
where the first term in the last line follows from~\eqref{E:deltav}.
Therefore
\[
|\nabla \cdot \nu - (\nabla_{X_i} \omega)(X_i) \cdot N|_{H^{3k-
\frac 52} (\p \Omega)} \le C |v|_{H^{3k-\frac 18} (\Omega)}.
\]
On the other hand, by decomposing $\nu$ into the tangential and
normal parts, one may calculate $\nabla \cdot \nu$ alternatively
\[
\nabla \cdot \nu = \CD\cdot \nu^\top + \kappa \nu^\perp + \nabla_N\,
\nu \cdot N= \CD\cdot \nu^\top + \kappa \nu^\perp + \nabla_N
(\nabla_{N_\CH}\, v \cdot N_{\CH}) - N \cdot D v(\CN (N))
\]
which along with the previous identity implies
\[
|\CD\cdot \nu^\top + \nabla_N (\nabla_{N_\CH}\, v \cdot N_{\CH}) -
(\nabla_{X_i} \omega)(X_i) \cdot N|_{H^{3k - \frac 52} (\p \Omega)}
\le C |v|_{H^{3k-\frac 18} (\Omega)}.
\]
Since $\omega$ is controlled by $E$ and  $\CD \cdot \nu^\top$
has been estimated in Step 2, we have
\[
|\nabla_N (\nabla_{N_\CH}\, v \cdot N_{\CH})|^2_{\dot H^{3k-\frac
52} (\p \Omega)} \le C(E + |v|_{H^{3k-\frac 18} (\Omega)}^2)
\]
The following decomposition trick on $\p \Omega$ has been used many
times in the basic estimates in Sections~\ref{pre} and~\ref{S:geo},
\begin{align*}
\nabla_N (\nabla_{N_\CH}\, v \cdot N_{\CH}) = &\CN (\nabla_N\, v
\cdot N) + \nabla_N \Delta^{-1} ( \nabla_{N_\CH} \Delta v \cdot
N_\CH + 2 \text{tr} ((D N_\CH)^* D v D N_\CH) \\
& + 2 D^2 v( N_\CH) \cdot D N_\CH + 2 D^2 v(\frac {\p N_\CH}{\p
x^i}, \frac \p{\p x^i}) \cdot N_\CH ).
\end{align*}
Using~\eqref{E:deltav} again, we have with a uniform $C>0$,
\[
|\CN(\nu^\perp)|^2_{H^{3k-\frac 52}(\p \Omega)} = |\CN (\nabla_N\, v
\cdot N)|^2_{\dot H^{3k-\frac 52}(\p \Omega)} \le C(E +
|v|_{H^{3k-\frac 18}(\Omega)}^2)
\]

From Step 3 and Step 4 above , we have
\[
|\nu|^2_{\dot H^{3k-\frac 32}(\p \Omega)} \le C (E + |v|_{H^{3k-
\frac 18}(\Omega)}^2)
\]
which implies the same estimate of the boundary data $\nabla_N v$ on
$\p \Omega$. Combining it with the Poisson
equation~\eqref{E:deltav}, we obtain
\[
|v|_{H^{3k} (\Omega)}^2 \le C (E + |v|_{H^{3k-\frac 18}(\Omega)}^2).
\]
Since $|v|_{H^{3k-\frac 18}(\Omega)} \le \beta |v|_{H^{3k}(\Omega)} + C_\beta |v|_{L^2(\Omega)}$,
proposition~\ref{P:energy} follows immediately. \hfill$\square$\\

In the proof of Theorem~\ref{T:STenergyE} and
Theorem~\ref{T:energyE}, we will need the following lemmas.

\begin{lemma} \label{L:Pi} For any $\Omega \in \Lambda_0$ with
$\kappa \in H^s (\p \Omega)$, $s\in [3k-\frac 52, 3k -\frac 32]$, we
have
\[
|\Pi|_{H^s (\p \Omega_t)} + |N|_{H^{s+1} (\p \Omega_t)} \le C(1+
|\kappa|_{H^s (\p \Omega_t)})
\]
for some $C>0$ uniform in $\Omega \in \Lambda_0$.
\end{lemma}

\begin{proof} We only need to prove the estimates for $\Pi$.
This is obvious if $n=2$. For $n\ge 3$, we use
identity~\eqref{E:DeltaPi}:
\[
- \Delta_{\p \Omega} \Pi = - \CD^2 \kappa + (|\Pi|^2 I- \kappa \Pi)
\Pi.
\]
Since $\Omega \in \Lambda_0$, we have
\[
|\Pi|_{H^{3k-\frac 52}(\p \Omega_t)} + |\kappa|_{H^{3k-\frac 52}(\p
\Omega_t)}\le C
\]
and it implies
\[
||\Pi|^2 I - \kappa \Pi|_{H^{s_1} (\p \Omega_t)} \le C, \qquad
\text{ where } \qquad s_1= \min \{6k-5-\frac{n-1}2, 3k-\frac 52\}
\]
as long as $3k-\frac52 \ne \frac{n-1}2$. Therefore, for $s_2\in
[3k-\frac 52, 3k -\frac 32]\backslash \{\frac{n-1}2\}$,
\[
|(|\Pi|^2 I- \kappa \Pi) \Pi|_{H^{s_3} (\p \Omega_t)} \le C
|\Pi|_{H^{s_2} (\p \Omega_t)} \qquad s_3 = \min\{ s_1, s_1+ s_2
-\frac{n-1}2\}.
\]
Since $s_3 +2 - s_2 \ge \min\{1, 3k-1-\frac n2\}>0$ and $s_1+2 > 3k
- \frac 32$, the estimate of $\Pi$ can be improved to $H^s(\p
\Omega_t)$ by bootstrap on $\Delta_{\p\Omega_t}^{-1}$. The
exceptional cases of the indices can be handled similarly.
\end{proof}

Although this regularity of $\Pi$ in terms of $\kappa$ follows
directly from Proposition~\ref{P:meank}, the point of this lemma is  that the
constant $C>0$ depends only on $\Lambda_0$, i.e.,  $\Omega\in H^{3k - \frac 12}$. This is also the point of
the following lemma.

\begin{lemma} \label{L:pressure}
For $\Omega \in \Lambda_0$ with $\p \Omega \in H^{3k}$ and $v \in
H^{3k}(\Omega)$, we have
\[
|\nabla_{N_\CH} p_{v,v}|_{H^{3k-\frac 12} (\Omega)} + |D^2
p_{v,v}|_{H^{3k-\frac 32} (\Omega)} \le C (1+ |\kappa|_{H^{3k-2} (\p
\Omega)}) |v|_{H^{3k} (\Omega)}^2.
\]
for some $C>0$ uniform in $\Omega \in \Lambda_0$.
\end{lemma}

\begin{proof} The idea of the proof is to use $p_{v,v}|_{\p \Omega_0}
=0$ and the identity $f =\CH(f|_{\p \Omega}) + \Delta^{-1} \Delta f$
for any $f : \Omega \to \RR$. On the one hand, notice on $\p
\Omega$,
\[
\nabla_N \nabla_{N_\CH} p_{v,v} = {\CN(N)} \cdot \nabla p_{v,v} +
D^2 p_{v,v} (N, N) =\CN(N)\cdot \nabla p_{v,v} + \Delta p_{v,v} -
\kappa \nabla_N p_{v,v},
\]
which, along Lemma~\ref{L:Pi} and~\eqref{E:pressure1}, implies
\[
|\nabla_N \nabla_{N_\CH} p_{v,v}|_{H^{3k-2} (\p \Omega)} \le C (1+
|\kappa|_{H^{3k-2} (\p \Omega)}) |v|_{H^{3k} (\Omega)}^2.
\]
On the other hand, using~Lemma \ref{L:Pi} and~\eqref{E:pressure1} as
well,
\begin{align*}
|\Delta \nabla_{N_\CH} p_{v,v}|_{H^{3k-\frac 52} (\Omega)} = &|-
N_\CH \cdot \nabla \text{tr} (Dv)^2 + 2 D^2 p_{v,v} \cdot D
N_\CH|_{H^{3k-\frac 52} (\Omega)}\\
\le& C (1+ |\kappa|_{H^{3k-2} (\p \Omega)}) |v|_{H^{3k} (\Omega)}^2.
\end{align*}
Therefore. we obtain the estimates on $\nabla_{N_\CH} p_{v,v}$

The estimate of $D^2 p_{v,v}$ is also achieved similarly. Firstly,
\[
|\Delta D^2 p_{v,v}|_{H^{3k-\frac 72} (\Omega)} = |D^2 \text{tr}
(Dv)^2|_{H^{3k-\frac 72} (\Omega)} \le C |v|_{H^{3k} (\Omega)}.
\]
For the boundary value of $D^2p_{v,v}$, we first consider
$D^2p_{v,v}(X, X)$ at $x \in \p \Omega$ with $X\in T_x \p \Omega$.
As usual, extend $X$ to a vector fields in a neighborhood of $x$ on
$\p \Omega$ by parallel transporting $X$ along geodesics emitting
from $x$. Thus,
\[
D^2 p_{v,v} (X, X) = \nabla_X \nabla_X p_{v,v} - \nabla_{\nabla_X X}
p_{v,v} = \nabla_N p_{v,v} \Pi(X, X).
\]
Also, we have
\begin{align*}
&D^2 p_{v, v} (N, N) = \Delta p_{v,v} - \kappa \nabla_N p_{v,v} = -
\text{tr}(Dv)^2 - \kappa \nabla_N p_{v,v}\\
&D^2 p_{v,v} (N, X) = \nabla_X \nabla_N p_{v,v} - \nabla_{\nabla_X
N} p_{v,v} = \nabla_X \nabla_N p_{v,v}.
\end{align*}
Therefore, from the above estimate on $\nabla_{N_\CH} p_{v,v}$, we
obtain the estimate on $D^2 p_{v,v}$.
\end{proof}

The following lemma is the most important observation of this paper,
which states that $J$ is a solution of the linearized Euler's
equation up to lower order terms. For the rest of this section, let
$Q$ denote a generic positive polynomial in $|v|_{H^{3k}
(\Omega_t)}$, $|\kappa|_{H^{3k-2} (\p \Omega_t)}$, and $\ep
|\kappa|_{H^{3k-1} (\p \Omega_t)}$ with coefficients depending only
on the set $\Lambda_0$.

\begin{lemma} \label{L:J}
Suppose a solution of the Euler's equation is given by $\Omega_t \in
\Lambda_0$ with $\p \Omega_t \in H^{s_0}$, $s_0 = 3k$ if $\ep=0$ or
$s_0=3k+1$ if $\ep>0$, and $v(t, \cdot) \in H^{3k} (\Omega_t)$, then
we have
\[
|\TD_t^2 J +\SR_0(v) J+ \ep^2 \CA J|_{H^{3k -3}(\Omega_t)} \le (1+
|\TD_t J|_{H^{3k -3}(\Omega_t)})Q.
\]
\end{lemma}

\begin{proof}
We start our proof by two simple observations. First we recall that  for any vector field $w$ defined on
$\Omega_t$ satisfying $\nabla \cdot w=0$, we write $p_{v,w} =-\Delta^{-1} \text{tr} (DvDw)$
\begin{equation} \label{E:IIE}
|p_{v, w}|_{H^{s+1} (\Omega_t)} \le C |v|_{H^{3k} (\Omega_t)}
|w|_{H^s (\Omega_t)}, \qquad
\quad s \in [0, 3k -1],
\end{equation}
where $C>0$ is uniform in $\Omega \in \Lambda_0$.
Second the first order derivative $\SD_tJ$ is lower order since by~\eqref{E:tdtJ} we have
\begin{equation} \label{E:tdtJ2}
|\TD_t J - \D_t J|_{H^s (\Omega_t)} \le Q \qquad  \qquad s\in [0, 3k -1],
\end{equation}
and by~\eqref{E:dtJ1} and~\eqref{E:dtk1}
\begin{equation}\label{E:dtJ3}
|\D_t J|_{H^{3k-3} (\p \Omega_t)} \le C |v|_{H^{3k} (\Omega_t)}.
\end{equation}

To verify the lemma we consider the expression for $\SD_{t}^2 J$ given in~\eqref{E:tdt2J1} and keep the least regular terms. Thus
\[
|\TD_t^2 J - \D_t^2 J|_{H^{3k -3} (\Omega_t)} \le Q + |\D_t \nabla p_{v, J}|_{H^{3k -3} (\Omega_t)}
\]
$\D_t \nabla p_{v, J}$ is given by~\eqref{E:tdt2J2} and can be estimated using~\eqref{E:IIE}, lemma~\ref{L:pressure} and Euler's equation
\[
|\D_t \nabla p_{v, J}|_{H^{3k -3} (\Omega_t)} \le Q + C |D\D_t v
\cdot DJ|_{H^{3k-4} (\Omega_t)} \le Q.
\]
To estimate  $\D^2_tJ$ given in~\eqref{E:dttJ1}, we use \eqref{E:pressure1} and Euler's equation to obtain
\[
|\D_t^2 J - \nabla \D_t^2 \kappa_\CH|_{H^{3k -3} (\Omega_t)}\le Q.
\]
To estimate $\nabla \D_{t}^2 \kappa_\CH$ we use ~\eqref{E:dt2chk} and~\eqref{E:dt2chk1}
\[
|\nabla \D_{t}^2 \kappa_\CH - \nabla \CH (\D_t^2 \kappa)|_{H^{3k -3}
(\Omega_t)} \le Q.
\]
Combine these inequalities, we obtain
\[
|\TD_t^2 J - \nabla \CH (\D_t^2 \kappa)|_{H^{3k -3} (\Omega_t)}\le
Q.
\]
The term $\D^2_{t} \kappa$ has been calculated explicitly
in~\eqref{E:dttk3} which yields,
\begin{equation} \label{E:tdt2J3}
|\TD_t^2 J + \nabla \CH (\Delta_{\p \Omega_t}\D_t v \cdot N - 2\ep^2
\Pi \cdot (D^\top|_{T\p\Omega_t} J))|_{H^{3k -3} (\Omega_t)}\le Q.
\end{equation}
To deal with $\Delta_{\p \Omega_t}\D_t v \cdot N $ we use Euler's equation to obtain that on the boundary
\[
\Delta_{\p \Omega_t}\D_t v \cdot N - 2\ep^2 \Pi \cdot
(D^\top|_{T\p\Omega_t} J))= - N \cdot\Delta_{\p \Omega_t} (\nabla
p_{v,v}) - \ep^2 \Delta_{\p \Omega} \CN (\kappa) +\ep^2 J \cdot
\Delta_{\p \Omega_t} N.
\]
The last term can be bounded by the identity
$
\Delta_{\p \Omega_t} N = -|\Pi|^2 N + \nabla^\top \kappa,
$
and Lemma~\ref{L:Pi}
$$
|\ep^2 J \cdot \Delta_{\p \Omega_t} N|_{H^{3k-\frac52} (\p\Omega_t)}
= |\ep^2 |\nabla^\top \kappa|^2- \ep^2 \CN(\kappa)
|\Pi|^2|_{H^{3k-\frac52} (\p\Omega_t)} < Q
$$
Substituting the above into~\eqref{E:tdt2J3} and Lemma~\ref{L:Pi}
\begin{equation} \label{E:tdt2J4}
|\TD_t^2 J - \nabla\CH(N \cdot\Delta_{\p \Omega_t} (\nabla p_{v,v}))
+ \ep^2 \CA J|_{H^{3k -3} (\Omega_t)} \le Q.
\end{equation}
As we are very close to the final desired form, the second term on
the above left side has to be related to $\SR_0$. Using
formula~\eqref{E:lapalace}, on $\p \Omega_t$, we have
\begin{align*}
&- N \cdot\Delta_{\p \Omega_t} \nabla p_{v,v} = -N \cdot \Delta
\nabla p_{v,v} + \kappa N \cdot \nabla_N \nabla p_{v,v} + N \cdot
D^2 (\nabla p_{v,v}) (N, N)\\
=& N \cdot \nabla (\text{tr}(Dv)^2) + \nabla_N (\kappa_\CH
\nabla_{N_\CH} p_{v,v} + D^2 p_{v,v} (N_\CH, N_\CH)) - \nabla_N
p_{v,v} J^\perp \\
&- \kappa \nabla p_{v,v} \cdot \nabla_N N_\CH -2 D^2 p_{v,v} (N,
\nabla_N N_\CH)
\end{align*}
Keeping the least regular terms and using lemma ~\ref{L:pressure}, implies
\begin{equation}\label{E:rrr}\begin{split}
|- N \cdot\Delta_{\p \Omega_t} \nabla p_{v,v} + \nabla_N p_{v,v}  J^\perp|_{H^{3k -\frac 52} (\p \Omega_t)}& \le Q \\
+  |\nabla_N & \lf(\kappa_\CH \nabla_{N_\CH} p_{v,v} + D^2 p_{v,v} (N_\CH, N_\CH)\rt) |_{H^{3k -\frac 52} (\p\Omega_t)}
\end{split}\end{equation}
Let $f=\kappa_\CH \nabla_{N_\CH} p_{v,v} + D^2 p_{v,v} (N_\CH,
N_\CH)$ defined on $\Omega_t$,  since $p_{v,v}|_{\p
\Omega_t}=0$ then
\[
|f|_{\p \Omega_t}|_{H^{3k -\frac32} (\p \Omega_t)} = |\Delta p_{v,v}
- \Delta_{\p \Omega_t} p_{v,v}|_{H^{3k -\frac32} (\p \Omega_t)} = |
\text{tr}(Dv)^2|_{H^{3k -\frac 32}(\p \Omega_t)} \le C |v|_{H^{3k}
(\Omega_t)}^2.
\]
Moreover it is easy to check from Lemma~\ref{L:pressure},
\[
|\Delta f|_{H^{3k-3} (\Omega_t)} \le Q.
\]
which implies  $|\nabla_N f|_{H^{3k-\frac 52} (\p\Omega_t)}  \le Q$.  Therefore~\eqref{E:rrr}   together with the definition of $\SR_0$ and the half derivative  behavior of $\nabla\CH$ implies
$$
\lf| \SR_0 J+ \nabla\CH(N \cdot\Delta_{\p \Omega_t} (\nabla p_{v,v}))\rt||_{H^{3k-3} (\Omega_t)} \le Q
$$
which together with~\eqref{E:tdt2J4} concludes the estimate in the statement of
the lemma.
\end{proof}

\noindent {\bf Proof of Theorems~\ref{T:STenergyE}
and~\ref{T:energyE}.} To prove Theorem~\ref{T:STenergyE},
in addition to Proposition~\ref{P:energy}, we need the following: a)
the estimates on the Lagrangian coordinates map and consequently
$\kappa\in{H^{3k-\frac 52} (\p \Omega_t)}$, b) estimates on $\omega
= D v - (D v)^*$, and  c) commutators involving $\D_t$, mostly have been done
in Section~\ref{S:geo}. In the following all constant $C>0$ will be determined only by the set
$\Lambda_0$. \\

\noindent {\it Estimate of the Lagrangian coordinate map $u(t, y)$.}
From our assumption on $v$, the ODE $u_t (t, y)= v(t, u(t, y))$
solving $u$ is well-posed. Since $u(t, \cdot): \Omega_0 \to
\Omega_t$ is volume preserving and $3k>\frac n2 +1$, it is easy to
derive, for any $s\in [0, 3k]$, and $f \in H^s (\Omega_t)$
\[
|f \circ u(t, \cdot)|_{H^s (\Omega_0)} \le C |f|_{H^s (\Omega_t)}
|u(t, \cdot)|_{H^{3k} (\Omega_0)}^s,
\]
where $C>0$ depends only on $s$. The proof follows simply from
induction and interpolation. By duality, for $s \in [0, 3k]$,
\[
|f \circ u(t, \cdot)|_{H^{-s} (\Omega_0)} \le C |f|_{H^{-s}
(\Omega_t)} |u(t, \cdot)^{-1}|_{H^{3k} (\Omega_0)}^s.
\]
Therefore,
\begin{equation}\label{E:LagE1}
|u(t, \cdot)-I|_{H^{3k} (\Omega_0)} \le C \int_0^t |v(t',
\cdot)|_{H^{3k} (\Omega_t)} |u(t', \cdot)|_{H^{3k} (\Omega_0)}^{3k}
\;dt',
\end{equation}
where $C>0$ depends only on $n$ and $k$. Let $\mu>0$ be a positive
large number to be specified later,
\begin{equation} \label{E:t0}
t_0 = \sup \{t \mid |v(t',\cdot)|_{H^{3k} (\Omega_t)} < \mu, \;
\forall t' \in [0, t]\},
\end{equation}
We have $t_0>0$ due to the continuity of $v(t, \cdot)$ in $H^{3k}
(\Omega_t)$. Then, for all $t \in [0, t_0]$,
\[
|u(t, \cdot)-I|_{H^{3k} (\Omega_0)} \le \mu \int_0^t |u(t',
\cdot)|_{H^{3k} (\Omega_0)}^{3k} \;dt'.
\]
Therefore, from ODE estimates, there exists $t_1>0$ and $C_2>0$
which depend only on $|v(0,\cdot)|_{H^{3k} (\Omega_t)}$ such that,
for all $0\le t \le \min\{t_0, t_1\}$,
\begin{equation}\label{E:LagE}
|u(t, \cdot)-I|_{H^{3k} (\Omega_0)} \le C_2 t.
\end{equation}
It implies the mean curvature estimate, for all $0\le t \le
\min\{t_0, t_1\}$,
\begin{equation}
|\kappa(t, \cdot)|_{H^{3k-\frac 52} (\p \Omega_t)} \le |\kappa(0,
\cdot)|_{H^{3k-\frac 52} (\p \Omega_0)} + C_3 t.
\end{equation}
Here it is easy to see from local coordinates constructed in
Section~\ref{pre} that $C_3$ is determined only by
$|v(0,\cdot)|_{H^{3k} (\Omega_t)}$ and the set $\Lambda_0$.
Therefore, there exists $t_2>0$ determined only by
$|v(0,\cdot)|_{H^{3k} (\Omega_t)}$ and the set $\Lambda_0$ such that
$\Omega_t \in \Lambda_0$ for $0\le t \le \min\{t_0, t_2\}$.\\

\noindent {\it Evolution of the curl $\omega = D v - (Dv)^*$.} From
equations (E) and~\eqref{E:dtnabla}, we have
\[
\D_t \omega = D \D_t v - (D \D_t v)^* + ((Dv)^*)^2 - (Dv)^2 =
((Dv)^*)^2 - (Dv)^2 = - (Dv)^* \omega - \omega Dv.
\]
It is clear how to obtain the estimate of $\omega$ in terms of $v$:
differentiating the above equation $3k-1$ times, multiplying it by
$D^{3k-1} \omega$ and integrating it on $\Omega_t$, we have
\begin{equation} \label{E:dtomega}
\frac d{dt} \int_{\Omega_t} |\omega|_{H^{3k-1} (\Omega_t)}^2 dx \le
C |v|_{H^{3k} (\Omega_t)} |\omega|_{H^{3k-1} (\Omega_t)}^2.
\end{equation}

\noindent {\it The commutator involving $\D_t$.} First,
from~\eqref{E:dtdeltat} and~\eqref{E:dtcn}, it is easy to verify
that, for any function $f$ defined on $\p \Omega_t$,
\begin{align}
& |[\D_t, \Delta_{\p\Omega_t}]|_{L(H^{s_1} (\p \Omega_t), H^{s_1-2}
(\p \Omega_t))} \le C |v|_{H^{3k} (\Omega_t)} \qquad &&s_1 \in
(\frac 72-3k, 3k -\frac12], \label{E:dtdeltat1}\\
& |[\D_t, \CN]|_{L(H^{s_2} (\p \Omega_t), H^{s_2-1}(\p \Omega_t))}
\le C |v|_{H^{3k} (\Omega_t)}  &&s_2 \in (1, 3k -\frac12]. \notag
\end{align}
To extend the range of $s_2$, we use the weak form of $[\D_t, \CN]$:
\[\begin{split}
\int_{\p \Omega_t} g[&\D_t, \CN]f dS = \int_{\p \Omega_t} g((D v) -
(D v)^*) (\nabla^\top f) \cdot N + g \CN(f)
\nabla_N v \cdot N \; dS \\
& + \int_{\Omega_t} g_\CH \nabla f_\CH \cdot \Delta v - 2 Dv (\nabla
f_\CH) \cdot \nabla g_\CH + \nabla g_\CH \cdot \nabla \Delta_D^{-1}
(2Dv \cdot D^2 f_\CH + \nabla f_\CH \cdot \Delta v) dx.
\end{split}\]
To conclude that the above estimate for $[\D_t, \CN]$ holds for
$s_2= \frac 12$. By interpolation,
\begin{equation} \label{E:dtcnW1}
|[\D_t, \CN]|_{L(H^s(\p \Omega_t), H^{s-1}(\p \Omega_t))} \le C
|v|_{H^{3k} (\Omega_t)}, \qquad s\in [\frac 12, 3k - \frac 12].
\end{equation}

\noindent {\it Evolution of $E$: first look.} Recall the expression
of $E_{RT}$ and $E$ written as $E=I_1+I_2$:
\begin{align*}
&I_1 =\frac {\ep^2}2 |\CA^{k -\frac12} J|_{L^2(\Omega_t)}^2 = \frac
{\ep^2}2 \int_{\p \Omega_t} \kappa \cdot \CN (- \Delta_{\p \Omega_t}
\CN)^{2k-1} \kappa \; dS \qquad I_2 =\frac 12 |\CA^{k -1} \SD_t
J|_{L^2(\Omega_t)}^2, \\
&E_{RT} =\frac 12 <(\SR_0(v)) \CA^{k -1} J, \CA^{k-1}
J>_{L^2(\Omega_t)} = \frac 12 \int_{\p \Omega_t} -\nabla_N p_{v,v}
|(-\CN \Delta_{\p \Omega_t})^{k-1} \CN (\kappa)|^2 \; dS.
\end{align*}
Also recall that
\[
\frac d{dt} dS = (\CD \cdot v^\top + \kappa v^\perp) dS.
\]
Since
\[
\nabla \cdot v|_{\p \Omega_t} = \CD \cdot v^\top + \kappa v^\perp +
\nabla_N v \cdot N=0
\]
then
\[
|\kappa v^\perp + \CD \cdot v^\top|_{H^{3k-\frac 32} (\p \Omega_t)}
= |\nabla_N v \cdot N|_{H^{3k-\frac 32} (\p \Omega_t)} \le C
|v|_{H^{3k} (\Omega_t)}
\]
and thus $\frac d{dt} dS$ would not complicate the estimates since
$3k\ge \frac n2 +\frac 32$.\\

\noindent {\bf I:} $\qquad |\frac d{dt} I_1 - \ep^2 <\CA^{k
-\frac12} J, \CA^{k -\frac12} \SD_t J>_{L^2(\Omega_t)}| \le Q.$ \\
To prove the inequality {\bf I}, we use~\eqref{E:dtdeltat1}
and~\eqref{E:dtcnW1} to obtain
\[
|\frac d{dt} I_1 - \ep^2 <(-\Delta_{\p\Omega_t}\CN)^{2k-1} \kappa,
\CN\D_t \kappa>_{L^2(\p\Omega_t)}| \le Q,
\]
and from~\eqref{E:IIE} and~\eqref{E:tdtJ2}, we have
\begin{equation} \label{E:tdtJ3}
|\SD_t J - \nabla \CH(\D_t \kappa)|_{H^{3k-\frac32} (\Omega_t)} \le
C |v|_{H^{3k} (\Omega_t)} |J| _{H^{3k -\frac32} (\Omega_t)}.
\end{equation}
It implies the estimate {\bf I} for $\frac d{dt} I_1$. \\

\noindent {\bf II:} $\qquad \frac d{dt} I_2 - <\CA^{k-1}\TD_t J,
\CA^{k-1}\TD_t^2 J>_{L^2(\Omega_t)}| \le Q$. \\
If $k=1$, which may happen when $n=2,3$,
\[
\frac d{dt} I_2 = <\TD_t J, \TD_t^2 J>_{L^2(\Omega_t)},
\]
where we used the fact $<\TD_t J, (\TD_t^2 -\D_t \TD_t)
J>_{L^2(\Omega_t)}=0$. If $k>1$,
\[
I_2 = \frac 12 \int_{\p \Omega_t} (\TD_t J)^\perp \cdot (-\Delta_{\p
\Omega_t}) (-\CN \Delta_{\p \Omega_t})^{2k-3} (\TD_t J)^\perp \; dS.
\]
From~\eqref{E:dtdeltat1} and~\eqref{E:dtcnW1}, we obtain
\[
|\frac d{dt} I_2 - \int_{\p \Omega_t} \D_t ((\TD_t J)^\perp) \cdot
(-\Delta_{\p \Omega_t}) (-\CN \Delta_{\p \Omega_t})^{2k-3} (\TD_t
J)^\perp \; dS| \le Q.
\]
On $\p \Omega_t$,
\[
\D_t (\TD_t J \cdot N) = (\D_t \TD_t J)\cdot N - \nabla_{(\TD_t
J)^\top} v \cdot N = (\TD_t^2 J)\cdot N + \nabla_N p_{v, \TD_t J} -
\nabla_{(\TD_t J)^\top} v \cdot N,
\]
which implies, along with~\eqref{E:IIE}, the estimate {\bf II} for
$I_2$.\\

\noindent {\bf III:}$\qquad |\frac d{dt} E_{RT} - <\CA^{k-1}
\SR_0(v) J, \CA^{k-1} \TD_t J>_{L^2 (\Omega_t)}|  \le Q$.\\
In general, for any function $f(t, \cdot)$ defined on $\p \Omega_t$
with $\Omega_t \in \Lambda_0$ satisfying $\int_{\p \Omega_t} f
dS=0$, we have
\[
\frac d{dt} \int_{\p \Omega_t} -\nabla_N p_{v,v} f^2 dS =
\int_{\p\Omega_t} -\nabla_N p_{v,v} f^2 (\kappa v^\perp + \CD \cdot
v^\top) - 2\nabla_N p_{v,v} f \D_t f- f^2 \D_t (\nabla_N p_{v,v}) dS
\]
Therefore, we obtain from~\eqref{E:pressure1},
\[
|\frac d{dt} \int_{\p \Omega_t} -\nabla_N p_{v,v} f^2 dS - \int_{\p
\Omega_t} - 2\nabla_N p_{v,v} f\D_t f - f^2 \D_t (\nabla_N p_{v,v})
\; dS| \le Q.
\]
Commuting $\D_t$ with $\nabla$ and $\Delta^{-1}$
by~\eqref{E:dtnabla} and~\eqref{E:dtdelta-1},
\begin{align*}
&\D_t (N \cdot \nabla p_{v,v}) = N \cdot (\nabla \D_t p_{v,v} -
(Dv)^* \nabla p_{v,v}) \\
=& - \nabla_N p_{v,v} \nabla_N v \cdot N + \nabla_N \Delta^{-1} (2Dv
\cdot D^2 p_{v,v} + \Delta v \cdot \nabla p_{v,v}) - \nabla_N
\Delta^{-1} (\D_t \text{tr} (Dv)^2),
\end{align*}
and using Euler's equation~\eqref{E:euler} to get
\[
\frac 12 \D_t \text{tr} (Dv)^2 = - \text{tr} (Dv)^3 - D^2 p_{v,v}
\cdot Dv - \ep^2 D^2 \kappa_\CH \cdot Dv.
\]
Therefore,
\[
|\frac d{dt} \int_{\p \Omega_t} -\nabla_N p_{v,v} f^2 dS - 2
\int_{\p \Omega_t} - \nabla_N p_{v,v} f\D_t f - \ep^2 f^2 \nabla_N
\Delta^{-1} (D^2 \kappa_\CH \cdot Dv) dS|\le Q
\]
Substituting $f = (-\CN \Delta_{\p \Omega_t})^{k-1} \CN (\kappa)$,
we obtain
\[
|\frac d{dt} E_{RT} - \int_{\p \Omega_t} - \nabla_N p_{v,v} (-\CN
\Delta_{\p \Omega_t})^{k-1} \CN (\kappa) \cdot \D_t (-\CN \Delta_{\p
\Omega_t})^{k-1} \CN (\kappa) \; dS| \le Q.
\]
From~\eqref{E:dtdeltat1}, \eqref{E:dtcnW1}, and~\eqref{E:tdtJ3},
\[
|\frac d{dt} E_{RT} - <\SR_0(v) \CA^{k-1} J, \CA^{k-1} \TD_t J>_{L^2
(\Omega_t)}| \le Q.
\]

In order to apply lemma~\ref{L:J} to the  estimate {\bf III}, we need to estimate
\[\begin{split}
I_3 \triangleq & <\CA^{k-1} \TD_t J, \CA^{k-1} \SR_0(v)
J>_{L^2(\Omega_t)} - <\SR_0(v) \CA^{k-1} J, \CA^{k-1}
\TD_t J>_{L^2(\Omega_t)}\\
=& \int_{\p \Omega_t} (\TD_t J)^\perp \cdot (-\Delta_{\p \Omega_t}
\CN)^{2k-2} (-\nabla_N p_{v,v} \CN\kappa) \\
&- (\TD_t J)^\perp \cdot (-\Delta_{\p \Omega_t} \CN)^{k-1}
\lf((-\nabla_N p_{v,v}) \CN (-\Delta_{\p \Omega_t} \CN)^{k-1}
\kappa\rt) \; dS.
\end{split}\]
Our strategy will be to move $\CN$ and the multiplication operator
by $-\nabla_N p_{v,v}$ in the first integrand by commuting them with
$\CN$ and $\Delta_{\p \Omega_t}$. Thus, we need to estimate
$[\Delta_{\p \Omega_t}, \nabla_N p_{v,v}]$, $[\CN, \nabla_N
p_{v,v}]$, and $[\Delta_{\p \Omega_t}, \CN]$. Using
Lemma~\ref{L:pressure}, \eqref{E:productN}, and~\eqref{E:deltacn3},
we have,
\begin{align*}
&|[\Delta_{\p \Omega_t}, \nabla_N p_{v,v}]|_{L(H^s(\p \Omega_t),
H^{s-\frac 32} (\p \Omega_t))} \le Q \qquad && s \in (3-3k,
\frac{n-1}2)\\
&|[\CN, \nabla_N p_{v,v}]|_{L(H^s(\p \Omega_t), H^{s-\frac 12} (\p
\Omega_t))} \le Q  && s \in [\frac 12,
\frac{n-1}2)\\
&|[\Delta_{\p \Omega_t}, \CN]|_{L(H^s(\p \Omega_t), H^{s-\frac 52}
(\p \Omega_t))} \le C && s \in (\frac 72-3k, 3k-1).
\end{align*}
Therefore we obtain $|I_3|\le Q$, which implies the estimate {\bf
III} on $\frac d{dt} E_{RT}$.\\

\noindent {\it Evolution of $\CE$.} Combining the estimates {\bf I,
II, III,} and Lemma~~\ref{L:J}, we obtain
\begin{equation} \label{E:energyE0}
|\frac d{dt}\CE| \le Q, \quad \text{ where }\quad Q=Q(|v|_{H^{3k}
(\Omega_t)}, |\kappa|_{H^{3k-2} (\p \Omega_t)},
\ep|\kappa|_{H^{3k-1} (\p \Omega_t)})
\end{equation}
is a polynomial with positive coefficients that depend  only on
$\Lambda_0$. This inequality on $[0,  \min\{t_0, t_2\}]$ where
 $t_0$ is defined in~\eqref{E:t0} and $t_2$ is determined only
by $|v(0, \cdot)|_{H^{3k}(\Omega_t)}$ and the set $\Lambda_0$. \\

\noindent {\it Proof of Theorem~\ref{T:energyE}.} Assume (RT) holds.
From Proposition~\ref{P:energy}, \eqref{E:energyE0}
implies~\eqref{E:energyE}. In addition, by choosing $\mu$ large
enough compared to the initial data, $t_0$ is bounded below by a
constant $t^*>0$ depending only on $|v(0,
\cdot)|_{H^{3k}(\Omega_0)}$ and the
set $\Lambda_0$. Theorem~\ref{T:energyE} follows immediately.\\

\noindent {\it Proof of Theorem~\ref{T:STenergyE}.} Assume $\ep>0$.
From inequality~\eqref{E:energyE0} and Proposition~\ref{P:energy},
we obtain
\begin{equation} \label{E:energy1}
E(t) -E(0) + E_{RT}(t) - E_{RT} (0) \le \int_0^t Q_\ep(|v|_{H^{3k}
(\Omega_{t'})}, |\kappa|_{H^{3k-1} (\p \Omega_{t'})})\; dt',
\end{equation}
where we use $Q_\ep$ to represent the dependence $Q$ on $\ep$. From
proposition~\ref{P:DeltaN} and ~\eqref{E:pressure1} we have
 $$
|E_{RT}| \le C |\nabla_N p_{v,v}|_{L^\infty (\p \Omega_t)}
|\kappa|_{H^{3k-2} (\p \Omega_t)}^2\le  C |v|^2_{H^{3k-\frac58} (\Omega_t)} |\kappa|_{H^{3k-2}(\p \Omega_t)}^2.
$$
Interpolating $v$ between $H^{3k}(\Omega_t)$ and $H^{3k-\frac
32}(\Omega_t)$ and $\kappa$ between $H^{3k-1}(\p \Omega_t)$ and
$H^{3k-\frac 52} (\p \Omega_t)$, we obtain from
Proposition~\ref{P:energy},
\[
|E_{RT}| \le \frac 12 E + C_1 (1+ |v|_{H^{3k-\frac32} (\Omega_t)}^m)
\]
for some integer $m>0$ where the constant $C_1$, which include
$|\kappa|_{H^{3k -\frac 52}(\p \Omega_t)}$, is determined only by
$E_0$ and the set $\Lambda_0$. Since $\D_t v = -\nabla p_{v, v}
-\ep^2 J$ is controlled by $E$ in $H^{3k- \frac 32} (\Omega_t)$ due
to Proposition~\ref{P:energy}, we can use the Lagrangian coordinate
map $u(t, \cdot)$ to estimate $|v(t, \cdot)|_{H^{3k -\frac 32}
(\Omega_t)} - |v(0, \cdot)|_{H^{3k -\frac 32} (\Omega_0)}$. Through
a similar procedure of the derivation of~\eqref{E:LagE} and using
Proposition~\ref{P:energy}, there exists $t_3>0$, depending only on
$|v(0, \cdot)|_{H^{3k}(\Omega_t)}$ and the set $\Lambda_0$ so that
for $0 \le t \le \min\{t_0, t_3\}$,
\[
\lf| |v(t, \cdot)|_{H^{3k -\frac 32} (\Omega_t)}^m - |v(0,
\cdot)|_{H^{3k -\frac 32} (\Omega_0)}^m \rt| \le \int_0^t Q_\ep \,
dt'
\]
for some  polynomial  $Q_\ep$ with positive coefficients.
Therefore,
\[
E_{RT} \le \frac 12 E + C_1 (1+ |v(0, \cdot)|_{H^{3k -\frac 32}
(\Omega_0)}^m) + \int_0^t Q_\ep  \, dt' \le  \frac 12 E +  C_1 +
\int_0^t Q_\ep \, dt',
\]
where $C_1$ is determined only by $|v(0, \cdot)|_{H^{\frac32 k
-\frac32}(\Omega_0)}$ and the set $\Lambda_0$. Thus
\[
E(\Omega_t, v(t, \cdot)) \le 2 E(\Omega_0, v(0, \cdot)) + C_\ep +
\int_0^t Q_\ep  \, dt'.
\]
By inserting the above inequality into~\eqref{E:energy1} and using
proposition~\ref{P:energy}, we obtain~\eqref{E:STenergyE}. By
choosing $\mu$ large enough compared to the initial data,
Theorem~\ref{T:STenergyE} follows. \hfill $\square$

\section{Examples of Lagrangian coordinate maps less smooth than
$\p \Omega_t$} \label{S:counterE}

In Section~\ref{S:estiE}, we established a priori estimates of the
free boundary Euler's equation. In particular, the estimates
indicate that $\p \Omega_t$ is $\frac 12$, if $\ep=0$, or $\frac
32$, if $\ep>0$, derivative smoother than the velocity fields
$v|_{\p \Omega_t}$. This is an improvement compared with the
regularity directly given by the ODE defining the Lagrangian
coordinate maps. It is very natural to guess that the Lagrangian
coordinate maps might be smoother as well. However, the following
examples show that the Lagrangian coordinate maps is only as smooth
as the velocity fields.

\subsection{Case 1: with surface tension} This is a relatively
easy case for the construction of the example since we need not to
worry about the sign assumption (RT): $-\nabla_N\, p_{v,v} >0$. The
example is given for $n=2$ and $\Omega_t = B_2 (1)$, the
2-dimensional open unit ball. In the polar coordinate, it is easy to
verify that
\[
v(t, r, \theta) = \Theta(r)\frac \p{\p \theta}, \quad p(t, r,
\theta) =\int_r^1 r' \Theta(r')^2 dr',   \qquad \text{supp}(\Theta)
\subset \subset B_2(1) \backslash \{0\}
\]
is a stationary solution of (E) and (BC) for $\ep=1$. The Lagrangian
coordinate map
\[
u(t, r_0, \theta_0) = (r_0, \theta_0 + t \Theta(r_0))
\]
is only as smooth as $v$.

\subsection{Case 2: without surface tension.} We will construct
an example in $\RR^2$ again, which satisfies the sign condition
$-\nabla_N p >0$. Consider the domain and the vector field in the
form
\begin{align*}
& \Omega_t = \{ (r, \theta) \mid r_1(t) < r < r_2(t)\},\quad r_1' =
\frac {a_1(t)}{r_1}, \quad r_2' = \frac {a_1(t)}{r_2}\\
& v(t, r, \theta) = \frac {a_1(t)}r \frac \p{\p r} + \Theta(t, r)
\frac \p{\p \theta},
\end{align*}
with the functions $a_1(t)$ and $\Theta(t, r)$ to be determined. In
the polar coordinate
\[
\nabla_{\frac \p {\p r}} {\frac \p {\p r}} =0, \quad \nabla_{\frac
\p {\p r}} \frac \p{\p \theta} = \nabla_{\frac \p {\p \theta}} \frac
\p {\p r} = \frac 1r \frac \p {\p \theta}, \quad \nabla_{\frac \p
{\p \theta}} \frac \p {\p \theta} = - r \frac \p {\p r}.
\]
One may calculate
\[
D v (\frac \p {\p r})  = - \frac {a_1(t)}{r^2} \frac \p {\p r} +
(\Theta_r + \frac \Theta r) \frac \p{\p \theta}, \quad Dv(\frac
\p{\p \theta}) = -r \Theta \frac \p {\p r} + \frac {a_1(t)}{r^2}
\frac \p{\p \theta}.
\]
Thus, it is clear that the above given form ensures $v$ is
divergence free:
\[
\nabla \cdot v = Dv(\frac \p {\p r}) \cdot \frac \p {\p r}  +
Dv(\frac 1r \frac \p {\p \theta}) \cdot \frac 1r \frac \p {\p
\theta} = 0.
\]
Moreover, the above calculation implies
\begin{align*}
\D_t v = &v_t + \nabla_v\, v = (\frac {a_1'(t)}r - \frac
{a_1(t)^2}{r^3} - r\Theta^2) \frac \p {\p r} + (\Theta_t + \frac
{a_1(t)}r \Theta_r + \frac {2a_1(t)}{r^2} \Theta)
\frac \p {\p \theta} \\
- \Delta p =& \text{tr} ((\nabla v)^2) = \frac {2a_1(t)^2}{r^4} - 2r
\Theta \Theta_r - 2\Theta^2 =\frac {2a_1(t)^2}{r^4} - \frac 1r (r^2
\Theta^2)_r.
\end{align*}
Therefore, $p$ must be radially symmetric, i.e. $p=p(t, r)$. It is
straight forward to calculate
\[
- \nabla p = - p_r \frac \p {\p r} = - (\frac {a_1(t)^2} {r^3} + r
\Theta^2 +\frac {a_2(t)}r) \frac \p {\p r}
\]
for some function $a_2(t)$. From the boundary condition $p(t,
r_1(t))=p(t, r_2(t))$, we need for the Euler's equation,
\begin{equation} \label{E:example1} \begin{cases}
a_1' = - a_2 = (\log \frac {r_2}{r_1})^{-1} \lf( \frac {a_1^2}2
(\frac 1{r_1^2} - \frac 1{r_2^2}) + \int_{r_1}^{r_2} r \Theta^2
dr\rt)\\
\D_t \Theta + \frac {2a_1(t)}{r^2} \Theta= \Theta_t + \frac
{a_1(t)}r \Theta_r + \frac {2a_1(t)}{r^2} \Theta =0.
\end{cases}\end{equation}
Let the Lagrangian coordinate map be $r = r(t, r_0, \theta_0)$ and
$\theta =\theta (t, r_0, \theta_0)$, with $r_0 \in [r_1(0),
r_2(0)]$. Due to the symmetry of the vector field, we have $r= r(t,
r_0)$ and $\theta = \theta_0 + \theta_1 (t, r_0)$ which satisfy
\[
\p_t r = \frac {a_1(t)}r \qquad \p_t \theta_1 = \Theta(t, r).
\]
It is easy to see
\[
r^2 = r_0^2 + 2A(t), \qquad A(t) = \int_0^t a_1(t') dt'.
\]
Thus, the system~\eqref{E:example1} is equivalent to
\begin{equation} \label{E:example2} \begin{cases}
A'' = (\log \frac {r_2(0)^2 +2A}{r_1(0)^2 + 2A})^{-\frac 12} \lf(
\frac {(A')^2}2 (\frac 1{r_1(0)^2 + 2A} - \frac 1{r_2(0)^2+2A}) +
\int_{r_1(0)}^{r_2(0)} r_0 (\p_t \theta_1(t, r_0))^2 dr_0\rt)\\
\p_{tt} \theta_1 (t, r_0) = - \frac {2A'}{r_0^2+ 2A} \p_t \theta_1 =
- (\p_t \log (r_0^2 + 2A)) \p_t \theta_1.
\end{cases}\end{equation}
The system~\eqref{E:example2} can be viewed as an ODE system for
$(A, \theta_1 (\cdot))$ on an open set in the Banach space $X = \RR
\times C^0([r_1(0), r_2(0)])$. Therefore the unique existence of
solution to this system is guaranteed. One can solve for
\[
\theta (t, r_0) = \theta_0 + \Theta(0, r_0) \int_0^t \frac
{r_0^2}{r_0^2+ 2A(\tau)} dt' \qquad \Theta (t, r) = \frac {r^2 -
2A}{r^2} \Theta(0, \sqrt{ r^2 -2A}).
\]
Thus, it is clear that the Lagrangian coordinate map is not smoother
than the velocity field. Finally, one notices that small $\Theta$
would ensure the sign condition $p_r(t, r_2(t)) <0$ and $p_r(t,
r_1(t))>0$ since small $\Theta$ would make $v = \frac {a_1(t)}r
\frac \p {\p r} + \Theta (t, r) \frac \p {\p \theta}$ a small
perturbation of the irrotational solution $v=\frac {\tilde a_1(t)}r
\frac \p {\p r}$, with a slightly different $\tilde a_1$, which
satisfies the sign condition.

\section{Appendix I: Basic estimates} \label{pre}

In free boundary problems, it often happens that the moving domain
$\Omega_t$ is of class $H^s$ and moves with an $H^{s_0}$ velocity
field with $s_0 \le s$. Moreover, the estimates usually involve
functions and vector fields defined on $\Omega_t$ and $\p \Omega_t$.
Therefore, in this section, we consider collections $\Lambda$ of
domains $\Omega$ which are $H^{s_0}$ close to some reference domain
and bounded in the $H^s$ class in some sense to be defined
rigorously. We will outline some basic estimates on functions
defined on $\Omega$ and $\p \Omega$ and some related operators.
Through tedious derivation, these estimates will be guaranteed to be
uniform for all $\Omega \in \Lambda$.

\subsection{Sobolev norms}\label{S:sn}
Let $\Omega \subset \RR^n$ be a bounded connected domain, viewing
$H^s(\Omega)$, $s\ge 0$, as a quotient space of $H^s(\RR^n)$, define
the norm
\[
|g|_{H^s(\Omega)} = \inf \{ |G|_{H^s(\RR^n)} \; : \; G \in
H^s(\RR^n), \; G|_\Omega = g\}
\]
where $|\cdot|_{H^s(\RR^n)}$ is defined through the Fourier
transform.  As usual, for $s\ge 0$, $H^s_0(\Omega)$ represents the closure of
$C^\infty_0(\Omega)$ in $H^s(\Omega)$ and $H^{-s}(\Omega)$ is
isometric to $(H^s_0 (\Omega))^*$. It is  important to note that with this definition of $H^s$ norm
the constants in Sobolev embedding  ($H^s \to
L^p$ or $C^\alpha$) are independent of $\Omega$.  The relationship between this definition of $H^s$ norm and the standard definition will be explored later on page \pageref{equi}.

%Moreover, for a
%uniformly bounded family of domains $\Omega_\lambda$,
%$|\Delta^{-1}|_{L(H^{-1}(\Omega_\lambda), H^1_0(\Omega_\lambda))}$
%is bounded uniform in $\Omega_\lambda$, where $\Delta^{-1}$ is
%under the homogeneous Dirichlet boundary condition.\\

\noindent {\bf $C^1 \cap H^2$ boundary $\p \Omega$.} To consider
functions defined on $\p \Omega$, let $\Omega \subset \RR^n$ be a
bounded connected domain with $\p \Omega$ of class $H^2 \cap C^1$.
Consider the local graph coordinates of $\p \Omega$ in orthonormal
frames. When two coordinate charts of this type overlap, it is easy
to verify that the transition map between these two local coordinate
maps is also of $C^1\cap H^2$. Therefore, on $\p \Omega$, the
definitions of spaces $C^1(\p \Omega) \cap H^2(\p \Omega)$ of scalar
functions and $C^0(\p \Omega) \cap H^1(\p \Omega)$ of $(k, l)$-type
tensors, though defined in local coordinates, are independent of the
choice of local coordinates. The Christofell symbols and the usual
geometric quantities of the hypersurface $\p \Omega$, such as the
second fundamental form and mean curvature are well-defined in
$L^2(\p \Omega)$ and the sectional curvature is in $L^1(\p \Omega)$,
for it is like the square of the second fundamental form. As these
will be referred to later, we give the explicit formula in local
coordinates here. Let $\{e_1, \ldots, e_n\}$ be an orthonormal frame
and $(x^1, \ldots, x^n)$ be the coordinates associated with this
frame. Suppose $\Omega$ locally is given by $x^n > f(x^1, \ldots,
x^{n-1})$ with $f \in H^2 \cap C^1$, then using $(x^1, \ldots,
x^{n-1})$ as the local coordinates, we have
\begin{align}
%& \text{The unit outer normal vector } = \frac
%1{\sqrt{1+|\nabla f|^2}} (\nabla f, -1); \notag\\
%& \text{The Lebeasgue
%measure } dS = \sqrt{1+|\nabla f|^2}; \notag\\
& \text{The Christofell symbols } \Gamma_{ij}^k = \frac {\p_k
f \p_{ij} f} {\sqrt{1+|\nabla f|^2}};\notag \\
& \text{The second fundamental form } \Pi(\frac {\p}{\p x^i}) \cdot
\frac {\p}{\p x^j} = - \frac {\p_{ij} f} {\sqrt{1+|\nabla
f|^2}}\notag \\
& \text{Mean curvature } \kappa = - \p_j \lf( \frac{\p_j f}{\sqrt{1
+ |\nabla f|^2}} \rt) = - \frac {\Delta f} {( 1+ |\nabla
f|^2)^{\frac 12}} + \frac {\p_i f \p_j f \p_{ij} f} {( 1+ |\nabla
f|^2)^{\frac 32}}; \label{E:meanK}\\
&\text{Sectional curvature } \R(\frac \p {\p x^i}, \frac \p {\p
x^j}, \frac \p {\p x^i}, \frac \p {\p x^j}) = \frac {\p_{ii} f
\p_{jj} f - (\p_{ij} f)^2} {1+|\nabla f|^2} \notag \\
& \text{Beltrami-Lapalace } \Delta_{\p \Omega} \phi = \text{tr}
\CD^2 \phi =\CD \cdot \nabla^\top \phi= \frac 1{\sqrt{1+|\nabla
f|^2}} \p_i \left( g^{ij} \sqrt{1+|\nabla f|^2} \p_j \phi \rt);
\notag
\end{align}
where the matrix $(g^{ij}) = ( \delta_{ij} + \p_i f \p_j f) ^{-1}$.
For a $C^0 \cap H^1$ tensor $T$ of $(k, l)$-type, the covariant
derivatives $\CD T$ is a $(k, l+1)$-type tensor in $L^2$. For any
$(k, l)$ tensor $T_1$ and $(k, l+1)$ tensor $T_2$ in $C^0 \cap H^1$,
one may verify, possibly through smooth approximations of $\p
\Omega$,
\[
\int_{\p \Omega} (\CD T_1) \cdot T_2 \; dS = \int_{\p \Omega}
\text{tr} \lf( T_1 \cdot \CD T_2 (\cdot, \cdot)\rt)\; dS.
\]
where, on the above right side, $\CD T_2(\cdot, \cdot)$ denotes the
$(k, l)$-type tensor define by $\CD T_2 (X, Y) ( \cdots) = (\CD_X
T_2) (Y, \cdots)$.

From this identity, for any $L^2$ tensor $T$, one can define $\CD
T$, in the distribution sense, as in the dual space of $C^0\cap H^1$
tensors. It is straightforward to verify that, for any $(k, l)$-type
tensor $T$ in $C^0 \cap H^1$, we have
\begin{equation} \label{E:H1}
\int_{\p \Omega} T \cdot \Delta_{\p \Omega} T  dS =- \int_{\p
\Omega} |\CD T|^2 dS.
\end{equation}
If $T, \CD T \in C^0 \cap H^1$, we have
\begin{align}
&\int_{\p \Omega} |\Delta_{\p \Omega} T|^2 dS = \int_{\p \Omega}
|\CD^2 T|^2 + \lf[ \CD_{\R (X_{j_1}, X_{j_2}) X_{j_1}}\; T - \R
(X_{j_1}, X_{j_2})\; \CD_{X_{j_1}}\; T\rt]
\cdot \CD_{X_{j_2}}\;T \notag\\
&\qquad \qquad \qquad\qquad - \frac 12 |\R (X_{j_1}, X_{j_2})T|^2
dS. \label{E:H2}
\end{align}
Here $\{X_1, \ldots, X_{n-1}\}$ is any pointwise orthonormal frame
of $\p \Omega$ which always appears in the trace form resulting in
the independence of the corresponding quantities of the choice of
the frame. The curvature acts on $T$ in the usual sense
\[
\R(X, Y) T = \CD_Y \CD_X T - \CD_X \CD_Y T - \CD_{[Y, X]} T = \CD^2
T (Y, X) - \CD^2 T (X, Y).
\]
Though $\R(X, Y)T$ seems to contain derivatives of $T$, one may
calculate
\[
(\R(X, Y) T)(X_1, X_2, \ldots) = - T(\R(X, Y)X_1, X_2, \ldots) - T(
X_1, \R(X, Y) X_2, \ldots) - \ldots.
\]
So the dependence of $\R(X, Y)T$ on $X, Y, T$ is only pointwise and
$\R$ vanishes if $T$ is a scalar function.

Since $I -\Delta_{\p \Omega}$ is a positive self-adjoint operator on
$L^2(\p \Omega)$, for $\phi :\p \Omega \to \RR$ and $r\ge 0$, we
define the Sobolev norm $|\cdot|_{H^r(\p \Omega)}$ on the surface
$\p \Omega$ as
\[
|\phi|^2_{H^r(\p \Omega)} = \int_{\p \Omega} |(I - \Delta_{\p
\Omega})^{\frac r2} \phi|^2 dS;\qquad \qquad |\cdot|_{L^2(\p
\Omega)} = |\cdot|_{H^0(\p \Omega)}.
\]
As usual, for $r\ge 0$, $|\cdot|_{H^{-r}(\p \Omega)}$ coincides with
$|\cdot|_{H^r(\p \Omega)^*}$. One may note here, since the
Christofell symbols $\Gamma_{ij}^k$ are only in $L^2$, while
$|T|_{H^1(\p \Omega)} < \infty$ implies that $\CD T \in L^2(\p
\Omega)$ from~\eqref{E:H1}, it does not imply that $T$ is in $H^1$
in local coordinates, except when $T$ is of $(0, 0)$-type, i.e. a
scalar function. Similarly, from~\eqref{E:H2}, $|T|_{H^2(\p \Omega)}
< \infty$ does not imply $\CD^2 T \in L^2 (\p \Omega)$.

\begin{remark}
When $n=2$, $\p \Omega$ is 1-dimension and $\Delta_{\p \Omega} =
\p_{ss}$ where $s$ is the arc length parameter, which is well
defined if $\p \Omega$ is in $W^{1,1}$. In fact, $\p \Omega \in H^s$,
$s > \frac 32$, is sufficient for the definitions of all the objects
intrinsic in $\p \Omega$.
\end{remark}

\noindent {\bf $H^s$ boundary $\p \Omega$, $s> \frac {n+1}2$.} For
the purpose of this paper, we assume $\p \Omega \subset \RR^n$ is in
$H^s$ with $s > \frac {n+1}2$. On the one hand, we defined the norm
$|\cdot|_{H^r(\p \Omega)}$ using the Beltrami-Lapalace $\Delta_{\p
\Omega}$ in the above. On the other hand, an obvious and traditional
way to define the Sobolev space $H^r(\p \Omega)$, $-s \le r\le s$
for scalar valued functions and $1-s \le r \le s-1$ for tensor
valued functions, is through local coordinate coverings of $\p
\Omega$ and the definition of the Sobolev space $H^r(\RR^{n-1})$.
From standard Sobolev inequalities, it is easy to see that the
latter definition of the spaces $H^r(\p \Omega)$ is actually
independent of local coordinates and naturally induces a topology on
$H^r(\p \Omega)$. In particular, when $r\ge 0$ is an integer,
straightforward calculation also shows that a function $f$ (or
tensor field) belongs to $H^r (\p \Omega)$ if and only if $f, \CD^r
f \in L^2 (\p \Omega)$. In fact, we have
\begin{proposition}\label{P4.1}
For $r \in [-s, s]$ ($r \in [1-s, s-1]$ for tensors), the norm
$|\cdot|_{H^r (\p \Omega)}$ is equivalent to the norm on $H^r(\p
\Omega)$ defined by using local coordinates.
\end{proposition}

The proof of this proposition follows from the standard elliptic
estimates using the local coordinates along with interpolation. In
particular, when $s$ is an integer, one may also prove it
geometrically. In fact, the proposition clearly holds for $r=1$ and
$r=2$ due to~\eqref{E:H1} and~\eqref{E:H2} and Sobolev inequalities.
When $r$ is an integer and $r\in [3, s]$ ($r \in [3, s-1]$ for
tensors), the proposition can be proved by using the following
identity
\begin{equation} \label{E:deltaD}
([\Delta_{\p \Omega} , \CD] T) (X) = \R (X, X_j) \CD_{X_j} \; T +
(\CD_{X_j} \R) (X, X_j) T + (\R(X, X_j) \CD \; T)(X_j).
\end{equation}
Finally, for non-integer or negative $r$, the proposition follows
from interpolation and duality. Another implication
of~\eqref{E:deltaD} is that $\CD: H^r(\p \Omega) \to H^{r-1} (\p
\Omega)$ is bounded and $|\cdot|_{L(H^r(\p \Omega), H^{r-1} (\p
\Omega))}$ depends only on $\R$ and its derivatives.

It is well known that the regularity of $\p \Omega$ can be
determined from the regularity of its mean curvature $\kappa$.
\begin{proposition}\label{P4.2}
Let $\Omega \subset \RR^n$ be a domain such that $\p \Omega \in
H^{s_0}$, $s_0 > \frac {n+1}2$. Suppose $|\kappa|_{H^{s-2} (\p
\Omega)} < \infty$ with $s > s_0$, then $\p \Omega \in H^s$.
\end{proposition}

Proposition~\ref{P4.2} can be proved by using local coordinates and
standard quasilinear estimates. Another proof can be based on the
following identity which is also used in a priori estimates.
Intuitively, let $e: \p \Omega \to \RR^n$ be the imbedding, then
$\Pi = -N \cdot \CD^2 e$ and $\kappa = - N \cdot \Delta_{\p
\Omega}\; e$, where $\Pi$ is viewed as a symmetric quadratic form.
Thus it is expected that the difference $\Delta_{\p \Omega} \Pi -
\CD^2 \kappa$ should be of lower order terms only. In fact,
\begin{equation}\label{E:DeltaPi}
- \Delta_{\p \Omega} \Pi = - \CD^2 \kappa + |\Pi|^2 \Pi - \kappa
\Pi^2.
\end{equation}
To prove the identity, at any $x \in \p \Omega$, let $\lambda_i$,
$i=1, \ldots, n-1$, be the eigenvalues of $\Pi(x)$ and $X_i$ be the
associated eigenvectors which form an orthonormal frame of $T_x \p
\Omega$. Parallel transport this frame to every base point in a
neighborhood of $x$ on $\p \Omega$ along the radial geodesics
emitting from $x$. From the construction, we have $\CD_{X_i} X_j
=0$, $[X_i, X_j]=0$, and $\Pi(X_i)= \lambda_i X_i$ at $x$. For any
$X =a^j X_j$ with constants $a^1, \ldots, a^{n-1}$, at $x$,
\begin{align*}
(\Delta_{\p \Omega} \Pi) (X, X)= & (\CD_{X_i} \CD_{X_i} \Pi) (X, X)
= \nabla_{X_i} ((\CD_{X_i} \Pi) (X, X)) = \nabla_{X_i} ((\CD_X \Pi)
(X_i, X))\\
= & (\CD_{X_i} \CD_X \Pi) (X_i, X) = (\CD_X \CD_{X_i} \Pi) (X, X_i)
+ (\R ( X, X_i) \Pi) (X_i, X)\\
=& \nabla_X ( (\CD_{X_i} \Pi) (X, X_i)) + \Pi ( \R (X_i, X) X_i, X)
+ \Pi (X_i, \R(X_i, X)X).
\end{align*}
For the first term at $x$, from the construction of our special
frame,
\[
\nabla_X ( (\CD_{X_i} \Pi) (X, X_i)) = \nabla_X ((\CD_X \Pi) (X_i,
X_i)) = \nabla_X \nabla_X \kappa - 2\Pi (X_i, \CD_X \CD_X X_i)=
\CD^2 \kappa (X, X).
\]
To calculate the remaining two terms, one may substitute $X= a^j
X_j$ and use $\Pi(X_j) = \lambda_j X_j$, the symmetry of $\R$, and
the following calculation
\begin{align*}
4\R(X_i, X_{j_1})X_i \cdot X_{j_2} =& \R(X_i, X_{j_1} + X_{j_2})X_i
\cdot (X_{j_1} +X_{j_2}) - \R(X_i, X_{j_1} -
X_{j_2})X_i \cdot (X_{j_1} -X_{j_2})\\
=& \Pi(X_i, X_i) \Pi (X_{j_1} + X_{j_2}, X_{j_1} + X_{j_2})-
\Pi(X_i, X_{j_1} + X_{j_2})^2  \\
&-\Pi(X_i, X_i) \Pi (X_{j_1} -X_{j_2}, X_{j_1} - X_{j_2})
+ \Pi(X_i, X_{j_1} - X_{j_2})^2 \\
=& 4\delta_{j_1j_2} \lambda_i \lambda_{j_1} - 4\delta_{ij_1}
\delta_{ij_2} \lambda_i^2.
\end{align*}
Equality~\eqref{E:DeltaPi} follows consequently. \\

\noindent {\bf $H^{s_0}$ neighborhoods of domains, $s_0 > \frac
{n+1}2$.} Given a domain $\Omega_*$ with $\p \Omega_*$ in $H^{s_0}$,
we will consider the set $\Lambda \triangleq \Lambda(\Omega_*, s_0,
s, L, \delta)$, $s\ge s_0$, of neighboring domains of $\Omega_*$,
given  in Definition~\ref{D:domainnbd}. From Proposition~\ref{P4.2}, every $\Omega \in \Lambda$ is in $H^s$.
Given $\Omega_*$ and sufficiently small $\delta>0$, in the
following, we will derive some estimates with bounds $C$ uniform in
$\Omega \in \Lambda$. Since $\p \Omega_*$ is compact, for any
$\sigma>0$, there exist $x_i \in \RR^n$ and $d, d_i \in (0,
\frac12]$, $i=1, \ldots, m$,
\begin{enumerate}
\item [(B1)] $B(\p \Omega_*, d) \subset \cup_{i=1}^m R_i(d_i)$
where each $R_i(\cdot)= \tilde R_i (\cdot) \times I_i (\cdot)
\subset \RR^n$ with $\tilde R_i(\cdot)$ and $I_i(\cdot)$ being an
open $(n-1)$-dimensional disk and an open perpendicular segment in
$\RR^n$, both centered at $x_i$ and of the given radius and half
length, respectively;
 \item [(B2)] For each $i$, $z= (z_1, \ldots, z_{n-1}, z_n)= (\tilde z,
 z_n)$
being an Euclidean coordinate system on $\tilde R_i (\cdot) \times
I_i (\cdot)$, there exists an $H^s$ function $f_{*i} : \tilde R_i
(2d_i) \to I_i$, so that
\begin{equation} \label{E:coodf}
|f_{*i}|_{C^0} < \sigma d_i, \qquad |D f_{*i}|_{C^0} < \sigma,
\qquad \text{ and} \qquad \Omega_* \cap R_i(2d_i)= \{ z^n > f_{*i}
(\tilde z) \}.
\end{equation}
\end{enumerate}

For any $\sigma >0$ with a fixed coordinate covering
$\{R_i(d_i)\}_{i=1}^m$ of $\p \Omega_*$ of the above type, it is
clear that, when $\delta >0$ is sufficiently small,
$\{R_i(d_i)\}_{i=1}^m$ is still a coordinate covering of any $\p
\Omega \in \Lambda$ satisfying (B1) and (B2) with coordinate
functions $\{f_i \in H^s\}_{i=1}^m$. This will provide us some
technical convenience in deriving estimates uniform in $\Omega \in
\Lambda$. The following proposition is a refinement of
Proposition~\ref{P4.2}.

\begin{proposition} \label{P:meank}
Given $\Omega_*$, there exists $\delta >0$ such that, for any $L>0$,
there exists $C>0$ such that the second fundamental form of any
$\Omega \in \Lambda$ satisfies
\[
|\Pi|_{H^{s-2}(\p \Omega)} < C.
\]
\end{proposition}

\begin{proof}
The proof follows simply from the standard elliptic estimates and we
will only give a sketch. With $\p \Omega \in H^s$ due to
Proposition~\ref{P4.2}, we will use the above coordinate covering
$\{R_i(d_i)\}_{i=1}^m$ and the coordinate functions $\{f_i \in H^s
(\tilde R_i(2d_i)) \}_{i=1}^m$, whose $H^{s_0} (\tilde R_i(2d_i))$
norms are uniformly bounded in $i$ and $\Omega$. Let $\gamma: [0,
+\infty) \to [0,1]$ be a smooth cut-off function supported on
$[0,\frac 32]$ and $\gamma|_{[0, \frac 54]} \equiv 1$. On each
$\tilde R_i(2d_i)$, let
\[
\gamma_i(\tilde z) = \gamma (\frac {|\tilde z|}{d_i}), \qquad
\kappa_i(\tilde z) = \gamma_i (\tilde z) \kappa(\tilde z, f_i(\tilde
z)), \qquad g_i = \gamma_i f_i,
\]
where $\kappa$ is the mean curvature of $\p \Omega$. It is clear
from the definition of $\Lambda$ that $|\kappa_i|_{H^{s_1-2} (\tilde
R_i (d_i))}$ is bounded uniformly in $i$ and $\Omega$ for $s_1 =
\min \{s_0 + 2, s\}$. From the mean curvature
formula~\eqref{E:meanK},
\[\begin{split}
- \Delta g_i + \frac {\p_{j_1} f_i \p_{j_2} f_i} {1+ |\nabla f_i|^2}
\p_{j_1 j_2} g_i = & ( 1+ |\nabla f_i|^2)^{\frac 12} \kappa_i -
\Delta \gamma_i f_i - 2 D \gamma_i \cdot D f_i \\
&+ \frac {\p_{j_1} f_i \p_{j_2} f_i} {1+ |\nabla f_i|^2} \lf(
\p_{j_1j_2} \gamma_i f_i + \p_{j_1} \gamma_i \p_{j_2} f_i + \p_{j_2}
\gamma_i \p_{j_1} f_i \rt)
\end{split}\]
Since $\gamma$ is supported on $[0, \frac 32]$, without loss of
generality, we may treat $f_i$ as compactly supported on the ball of
radius $\frac{7d_i}4$ because $f_i$ can always be replaced by
$\tilde \gamma (\frac {|\tilde z|}{d_i}) f_i(\tilde z)$ where
$\tilde \gamma$ is a cut-off function supported on $[0, \frac 74]$
and $\tilde \gamma|_{[0,\frac 32]} =1$. By partition of the unity
and the Inverse Function Theorem, $(f_1, \ldots, f_m)$ can be
expressed by $g = (g_1, \ldots, g_m)$ with the same regularity and
similar estimates. Thus, dividing both sides of the above equation
by $( 1+ |\nabla f_i|^2)^{\frac 12}$, it can be rewritten as
\[
- \Delta g + A_{j_1 j_2} (x, g, \p g) \p_{j_1 j_2} g = \kappa + G
(x, g, \p g),
\]
where $A_{j_1j_2}$, $G_1$, $G_2$ are smooth in their arguments and
$A_{j_1j_2} \le C \sigma^2$ so the left side is uniformly elliptic.
In this form, the estimate on $|g|_{H^{s_1} (\RR^{n-1})}$ uniform in
$\Omega$ is obtained following the standard theory of quasilinear
elliptic equations. If $s_1 < s$, this procedure can be carried out
again with $s_0$ replaced by $s_1 =s_0 +2$. Thus the desired uniform
estimates on $f_i$ in $H^s$ follow by repeating this procedure.
\end{proof}

\begin{remark}
1) One can also prove Proposition~\ref{P:meank} based
on~\eqref{E:DeltaPi}.\\
2) A more careful estimate can be found in Lemma~\ref{L:Pi}, when a
more detailed relationship between $|\Pi|_{H^{s-2} (\p \Omega)}$ and
$|\kappa|_{H^{s-2} (\p\Omega)}$ is given under certain conditions.
\end{remark}

Using~\eqref{E:H1}, \eqref{E:H2}, ~\eqref{E:deltaD}, and the above
uniform estimate on $\Pi$, which implies the uniform estimate on the
curvature $\R$, it is easy to prove that, for any tensor $T \in
H^r(\p \Omega)$, $r \in [2-s, s-1]$ ($r \in [1-s, s]$ for scalars),
we have
\begin{equation} \label{E:cdnorm}
|\CD T|_{H^{r-1} (\p \Omega)} \le C |T|_{H^r (\p \Omega)}
\end{equation}
for some $C$ uniform in $\Omega \in \Lambda$.

From the uniform estimates on those (uniformly fixed) local
coordinates derived in the above proof, it is also clear that the
constants in the Sobolev inequalities (e.g. $|\cdot|_{H^s(\p
\Omega)}$ to $L^p(\p \Omega)$ or $C^\alpha(\p \Omega)$ for $s \le
k$) are uniform in $\Omega \in \Lambda$. The two most used
inequalities in this paper are for $f\in H^{s_1} (\p \Omega)$ and $g
\in H^{s_2} (\p \Omega)$, $s_1 \le s_2$,
\begin{align*}
& |fg|_{H^{s_1 +s_2 -\frac {n-1}2} (\p\Omega)} \le C |f|_{H^{s_1}(\p
\Omega)} |g|_{H^{s_2}(\p \Omega)}, \quad &&\text{ if } s_2 < \frac
{n-1}2 \text{ and } 0< s_1+s_2 \\
& |fg|_{H^{s_1} (\p\Omega)} \le C |f|_{H^{s_1}(\p \Omega)}
|g|_{H^{s_2}(\p \Omega)}, && \text{ if } s_2 > \frac {n-1}2 \text{
and } 0\le  s_1+s_2.
\end{align*}
Similar inequalities hold for $f$ and $g$ defined in $\Omega$.

\subsection{Dirichlet-Neumann operator}

Given $\Omega_*$, in order to study the Dirichlet-Neumann operator
for domains $\Omega \in \Lambda\triangleq \Lambda(\Omega_*, s_0, s,
L, \delta)$, we need to first construct local coordinate maps on
each $R_i(2d_i)$ for each $\Omega$, which flatten $\p \Omega$ and
have estimates uniform in $\Omega \in \Lambda$, based on the above
coordinates functions of $\p \Omega$.\\

\noindent {\bf Local coordinates and partition of the unit.} From
Proposition~\ref{P:meank}, $\p \Omega \cap R_i(2d_i)$ is represented
as the graph of an $H^s$ function $f_i: \tilde R_i(2d_i) \to \RR$.
Let $\phi = \gamma_i f_i$ with $\gamma_i$ defined in the previous
proof. A standard way to extend $\phi$ to a function $\Phi \in H^{s
+\frac 12} (\RR^n)$ is through the Fourier transform with an
appropriate constant $a$:
\begin{equation}\label{E:ext}
\hat \Phi (\xi^1, \ldots, \xi^n) = a \frac {(1+ (\xi^1)^2 + \ldots +
(\xi^{n-1})^2)^s}{(1+ (\xi^1)^2 + \ldots + (\xi^n)^2)^{s + \frac
12}} \hat \phi (\xi^1, \ldots, \xi^{n-1})
\end{equation}
Since $s + \frac 12> \frac n2 + 1$ and $\Phi$ is bounded in $H^{s+
\frac 12}$ uniformly in $i$ and $\Omega$, $|D \Phi|_{C^0}$ is also
uniformly bounded. Therefore, there exists $b>0$ so that
\begin{equation} \label{E:coord1}
H_i(z^1,\ldots, z^n) = (z^1, \ldots, z^{n-1}, bz^n + \Phi(z^1,
\ldots, z^n))
\end{equation}
is a diffeomorphism so that $|H_i|_{H^{s+\frac 12}(\RR^n)}$ and
$|(H_i)^{-1}|_{H^{s+\frac 12}(\RR^n)}$ are bounded uniformly in $i$
and $\Omega$. Let $G_i=(H_i)^{-1}$ and $g_i(z)$ be the $n$-the
component of $G_i$, then $|g_i|_{H^{s+\frac 12}(\RR^n)}$, $\p_{z_n}
g_i$, and $(\p_{z_n} g_i)^{-1}$ are bounded uniformly in $i$ and
$\Omega$. Obviously, there exists a uniform $\delta_*>0$ so that
\[
( \tilde R_i(\frac 54 d_i) \times I_i(\frac 54 \delta_* d_i)) \cap
\Omega = \tilde R_i(\frac 54 d_i) \times I_i(\frac 54\delta_* d_i))
\cap \{ g_i >0\}.
\]

Based on the local coordinate maps, we can construct partition of
the unit satisfying estimates uniform in $\Omega$ if $\delta$ is
small. In fact, take $\gamma, \xi: C^\infty ([0, +\infty), [0,1])$
so that supp$(\gamma) \subset [0,\frac 54]$, $\gamma|_{[0, \frac
98]} \equiv 1$, $\xi(r) = r$ for $r \ge \frac 23$, and $\xi|_{[0,
\frac 13]}\equiv \frac 13$. Define
\[
\tilde \gamma_{*i} (z) = \gamma (\frac {|\tilde z|}{d_i}) \gamma (
\frac {|z_n|} {\delta_* d_i}), \quad \eta = \xi \circ \Sigma_{i=1}^m
(\tilde \gamma_{*i} \circ G_i), \quad \gamma_{*i} = \frac {\tilde
\gamma_{*i} \circ G_i} \eta, \quad \gamma_{*0} = (1 - \Sigma_{i=1}^m
\gamma_{*i})\chi(\Omega).
\]
It is straight forward to verify that $\gamma_{*0}, \gamma_{*1},
\ldots, \gamma_{*m} \in H^{s+\frac 12}(\RR^n, [0, 1])$ satisfy
\[
|\gamma_{*i}|_{H^{s+\frac 12}(\RR^n)} \le C, \qquad
\text{supp}(\gamma_{*i}) \subset \tilde R_i(\frac 54 d_i) \times
I_i(\frac 54 \delta_* d_i)),
\]
for $i=1, \ldots, m$, and $(\Sigma_{i=0}^m \gamma_{*i})|_\Omega
\equiv
1$.\\
\label{equi}
\begin{remark}
Using the above local coordinates and partition of unity we can establish the equivalence of the standard  $H^\ell$ norm and the norm  given in definition \ref{S:sn} for integer $\ell \in (-s -\frac 12, s + \frac 12)$.  The ratio of the two norms is bounded above and below by two constants depending only on $\Lambda$.
\end{remark}
\noindent {\bf Trace and Harmonic extension.} Let $s_1 \in (\frac
12, s + \frac 12]$. Using the partition of the unit and the above
local coordinates, it is straight forward to obtain the trace
operator estimate
\begin{equation} \label{E:trace}
| \; (\Psi|_{\p \Omega})\;|_{H^{s_1-\frac 12} (\p\Omega)} \le C
|\Psi|_{H^{s_1}(\Omega)}
\end{equation}
for any $\Psi \in H^{s_1} (\Omega)$ where $C>0$ is uniform in
$\Omega \in \Lambda$.

In order to obtain the estimate on the Harmonic extension operator,
we first construct an extension for convenience. Let $s_2 \in (0,
s]$ and $\psi\in H^{s_2} (\p \Omega)$. Take the same auxiliary
functions $\gamma$ and $\xi$ used above. For each $1\le i \le m$,
let $\phi_i (\tilde z) = \gamma (\frac {|\tilde z|}{d_i}) \psi( H_i
(\tilde z, 0))$ and $\Phi_i(z)$ be the extension of $\phi_i$,
defined in the way of~\eqref{E:ext}. Let
\[
\tilde \Phi_i (z) = \Phi_i(z) \gamma (\frac {|\tilde z|}{d_i})
\gamma ( \frac {|z_n|} {\delta_* d_i}), \qquad \qquad \Psi_1 =
\Sigma_{i=1}^m \tilde \Phi_i \circ G_i
\]
where $\Psi_1$ can be viewed as a function defined on $\RR^n \supset
\Omega$. Let
\[
\eta_i (z) = \gamma (\frac {|\tilde z|}{d_i})^2 \gamma ( \frac
{|z_n|} {\delta_* d_i}), \qquad \eta = \xi \circ \Sigma_{i=1}^m
(\eta_i \circ G_i), \qquad \Psi = \frac {\Psi_1}{\eta}.
\]
It is easy to verify that $\Psi \in H^{s_2+\frac 12} (\RR^n)$ is an
extension of $\psi \in H^{s_2} (\p \Omega)$ satisfying the estimate
\begin{equation} \label{E:ext1}
|\Psi|_{H^{s_2 +\frac 12} (\R^n)} \le C |\psi|_{H^{s_2} (\p \Omega)}
\end{equation}
with $C$ uniform in $\Omega \in \Lambda$.

Using the partition of the unit and the local coordinates we
constructed above and following the standard procedure, we have

\begin{lemma} \label{L:Delta}
There exists $C>0$ which depends only on the set $\Lambda$ so that,
for $s_1 \in [\frac 12, s]$
\[
|\Delta^{-1}|_{L(H^{s_1-\frac 32} (\Omega), H^{s_1+\frac 12}
(\Omega))} + |\CH|_{L(H^{s_1} (\Omega), H^{s_1+\frac 12} (\Omega))}
\le C.
\]
\end{lemma}

\noindent {\bf Dirichlet-Neumann operator.} Following from the above
estimate, the Dirichlet-Neumann operator $\CN: H^{s_1}(\p \Omega)
\to H^{s_1-1} (\p \Omega)$ can be defined and it has a uniform bound
for $s_1 \in (1, s]$. In fact, we can extend $\CN$ into a weaker
form defined on $H^{s_1} (\p \Omega)$ for $s_1\ge \frac 12$. Given
$f \in H^{\frac 12} (\p \Omega)$, define $\CN(f) \in H^{-\frac 12}
(\p \Omega)$ as
\[
<\psi, \CN(f)> = \int_\Omega \nabla f_\CH \cdot \nabla \psi_\CH dx
\]
for any $\psi \in H^{\frac 12} (\p \Omega)$. It is easy to prove
that
\begin{enumerate}
\item[1)] $\CN$ is self-adjoint in $L^2(\p \Omega)$ with compact resolvent;
 \item[2)] the kernel $\ker(\CN) = \{ \text{const} \}$;
 \item[3)] $ C |f|_{H^{\frac 12} (\p \Omega)} \ge
|\CN(f)|_{H^{-\frac 12}(\p \Omega)} \ge \frac 1C |f|_{H^{\frac 12}
(\p \Omega)}$ for any $f$ satisfying $\int_{\p \Omega} f dS =0$.
\end{enumerate}
The first inequality of 3) follows from the uniform bound on $\CH$.
In order to prove the second inequality in (3), one notices that
\[
|f|_{H^{\frac 12} (\p \Omega)} |\CN(f)|_{H^{-\frac 12} (\p \Omega)}
\ge |<f, \CN(f)>| =\int_\Omega |\nabla \CH(f)|^2 dx.
\]
From the estimate of the trace operator, we only need, for any $f$
satisfying $\int_{\p \Omega} f dS =0$,
\[
|\CH(f)|_{L^2(\Omega)} \le C |\nabla \CH(f)|_{L^2(\Omega)}
\]
with a constant $C$ uniform in $f$ and $\Omega$. This inequality can
be proved by a compactness argument. Thus, by duality and
interpolation, $\CN$ can be extended to $H^{s_1} (\p \Omega)$ for
all $s \in [1-s, s]$ and $|\CN|_{L(H^{s_1} (\p \Omega), H^{s_1-1}
(\p \Omega))}$ is bounded uniformly in $\Omega$. Moreover, for $f
\in H^s (\p \Omega)$ with $\int_{\p \Omega} f dS =0$, we can obtain
$|f|_{H^s(\p \Omega)} \le C |\CN(f) |_{H^{s-1} (\p \Omega)}$ with
$C$ uniform in $\Omega$. The proof is simply the elliptic estimate
under the Neumann boundary condition -- very much similar to  the
derivation of the harmonic extension estimate, except in the first
step, instead of using~\eqref{E:ext}, we need to construct $F$ with
$|F|_{H^{s+\frac 12} (\Omega)} \le C|\CN(f)|_{H^{s-1} (\p \Omega)}$
and $\nabla_N F = \CN(f)$ on $\p \Omega$, by using a slightly
different formula of the same fashion. Therefore, from
interpolation, we have, for any $s_1 \in [\frac 12, s]$,
\[
|f|_{H^{s_1}(\p \Omega)} \le C |\CN(f) |_{H^{s_1-1} (\p \Omega)},
\qquad \text{ if } \int_{\p \Omega} f dS =0
\]
with $C$ uniform in $\Omega$. moreover, this inequality holds for
$s_1\in [1-s, s]$ by duality. Based on these estimates, we can use
$I +\CN$ to define the Sobolev norms which are equivalent to those
defined by using $I - \Delta_{\p \Omega}$ uniformly in $\Omega$,
i.e.

\begin{proposition} \label{P:cn}
For $s_1 \in [-s, s]$, the norms on $H^{s_1}(\p \Omega)$ defined by
interpolating $I -\Delta^T$ and $I +\CN$ are equivalent, i.e.
\[
\frac 1C (I - \Delta_{\p \Omega})^{\frac {s_1}2} \le  (I +
\CN)^{s_1} \le C (I - \Delta_{\p \Omega})^{\frac {s_1}2}
\]
with $C$ uniform in $\Omega \in \Lambda$.
\end{proposition}

Furthermore, for $s_1 \in [-s, s-1]$,
\[
\CN^{-1}: \dot H^{s_1} (\p \Omega) \to \dot H^{s_1+1} (\p \Omega),
\qquad \dot H^{s_1} (\p \Omega) = \{f \in H^{s_1} (\p \Omega)\;
\mid\; \int_{\p \Omega} f dS =0 \}
\]
is well defined and bounded uniformly in $\Omega$. $\CN^{-1}$
defined on $\cdot H^{-\frac 12} (\p \Omega)$ induces the solvability
of the Lapalace equation with Neumann boundary data given in
$H^{-\frac 12} (\p \Omega)$.

To demonstrate that $\CN$ behaves like differentiation, we give the
following ``product rule". Given functions $f$ and $g$ defined on
$\p \Omega$. Since
\[
f_\CH g_\CH - \Delta^{-1} \Delta (f_\CH g_\CH) = \CH(f_\CH g_\CH
|_{\p \Omega}) = \CH (fg) \qquad \text{ in } \Omega,
\]
we obtain
\begin{equation} \label{E:productN}
\CN(fg) = f \CN(g) + g \CN(f) - 2\nabla_N \Delta^{-1} (\nabla f_\CH
\cdot \nabla g_\CH).
\end{equation}
Since $\CN$ is like differentiation, coordinate independent, and
self-adjoint, appearing naturally in the Euler's equation, it is
sometimes convenient to express the Sobolev norms on $\p \Omega$
by $\CN$.\\

\noindent {\bf Relationship between $\CN$ and $\Delta_{\p \Omega}$.}
In addition to just the comparison between the norms of $\Delta_{\p
\Omega}$ and $\CN$, we will prove that $\CN$ is simply equal to
$(-\Delta_{\p \Omega})^{\frac 12}$ plus lower order terms. This
improves the previous estimates and makes the estimates of some
Sobolev norms using $\CN$ more convenient. From the identity
\begin{equation}\label{E:lapalace}
\Delta \psi= \Delta_{\p \Omega} \psi + \kappa \nabla_N \psi + D^2
\psi (N, N)\qquad x \in \p \Omega
\end{equation}
for any smooth function $\psi$ on $\Omega$. Recall that $N_\CH(x)$
and $\kappa_\CH (x)$, $x \in \Omega$, denote the harmonic extension
of the unit outward normal vector and the mean curvature of $\p
\Omega$. Given smooth $f: \p \Omega \to \RR$, at any $x \in \p
\Omega$,
\begin{align*}
D^2 f_\CH (N, N) =& \nabla_N \nabla_{N_\CH} f_\CH - \nabla f_\CH
\cdot \nabla_N N_\CH \\
=& \nabla_N \lf(\CH((\nabla_N f_\CH)|_{\p \Omega}) + (-\Delta)^{-1}
(-\Delta) (\nabla_{N\CH} f_\CH)\rt)
- \CN(N) \cdot (\CN(f) N + \nabla^\top f)\\
=& \CN^2 (f) - 2\nabla_N (- \Delta)^{-1} (D N_\CH \cdot D^2 f_\CH) -
\CN(N) \cdot (\CN(f) N + \nabla^\top f)
\end{align*}
which implies
\begin{equation} \label{E:deltacn1}
(-\Delta_{\p \Omega} - \CN^2) f = \kappa \CN(f) - 2\nabla_N (-
\Delta)^{-1} (D N_\CH \cdot D^2 f_\CH) - \CN(N) \cdot (\CN(f) N +
\nabla^\top f).
\end{equation}

\begin{proposition} \label{P:DeltaN}1) For $s > \frac {n+3}2$, there
exists $C>0$ uniform in $\Omega \in \Lambda$ such that we have
\[
|\Delta_{\p \Omega} + \CN^2|_{L(H^{s'} (\p \Omega), H^{s'-1} (\p
\Omega))} \le C, \qquad s' \in [2-s, s-1].
\]
2) For $s \in (\frac {n+1}2, \frac {n+3}2)$ and $s>2$, there exists
$C>0$ uniform in $\Omega \in \Lambda$ such that we have
\[
|\Delta_{\p \Omega} + \CN^2|_{L(H^{s'} (\p \Omega),
H^{s'-\frac{n+5}2 +s} (\p \Omega))} \le C, \qquad s' \in (2-s, \frac
{n+1}2).
\]
\end{proposition}

\begin{proof}
For $s'> \frac {n+5}2 -s$, the above inequalities follow directly
from~\eqref{E:deltacn1} and the estimates on $\CH$ and $\CN$. Thus,
by duality and interpolation, we only need to consider $s' = \frac
12$ or $s'=\frac \alpha2$ in each case, respectively. Let $f, g: \p
\Omega \to \RR$ be smooth and harmonically extend into $\Omega$.
Equality~\eqref{E:deltacn1} yields
\[
\int_{\p \Omega} g (-\Delta_{\p \Omega} - \CN^2) f dS = \int_{\p
\Omega} \kappa g \CN(f) + g \CN(N) \cdot (\CN(f) N + \nabla^\top f)
dS  - 2 \int_\Omega D N_\CH (\nabla g_\CH) \cdot \nabla f_\CH dx,
\]
which is sufficient to establish the estimate.
\end{proof}

\begin{coro}
The proposition implies the commutator estimates
\begin{equation}\label{E:deltacn2}
|[\Delta_{\p \Omega}, \CN]|_{L((H^{s'} (\p \Omega), H^{s'-2} (\p
\Omega))} \le C, \qquad s' \in [3-s, s-1]
\end{equation}
if $s>\frac {n+3}2$ and
\begin{equation}\label{E:deltacn3}
|[\Delta_{\p \Omega}, \CN]|_{L((H^{s'} (\p \Omega),
H^{s'-\frac{n+7}2 +s} (\p \Omega))} \le C, \qquad s' \in (3-s, \frac
{n+1}2),
\end{equation}
if $s \in (\frac {n+1}2, \frac {n+3}2)$ and $s>2$.
\end{coro}

We need the following abstract result for a more careful estimate on
$\CN$.

\begin{proposition}\label{P:squareroot}
Let $X$ be a Hilbert space and $A$ and $B$ be (possibly unbounded)
self-adjoint positive operators on $X$ so that $A^{-1}B$ and
$AB^{-1}$ are bounded. Suppose $K=A^2 - B^2 $ satisfies that
$KB^{-\alpha}$ is bounded with $\alpha \in [0, 2)$, then $(A-B)
B^{1-\alpha}$ is bounded.
\end{proposition}

\begin{proof}
Let $R= A-B$. Calculating $(B+R)^2 = B^2 +K$, we obtain
\[
-BR - RB = R^2 -K,
\]
which implies
\[
\frac d{dt} (e^{-Bt} R e^{-Bt}) = e^{-Bt} (R^2 -K) e^{-Bt} \ge -
e^{-Bt} K e^{-Bt} \ge - C e^{-Bt} B^\alpha e^{-Bt}.
\]
Therefore,
\[
R \le C_1\int_0^\infty e^{-Bt} B^\alpha e^{-Bt} dt = \frac {C_1}2
B^{\alpha -1}
\]
Calculating $A^2 = (A-R)^2 +K$ with a similar procedure, we obtain
\[
R \ge - C_2\int_0^\infty e^{-At} A^\alpha e^{-At} dt = - \frac
{C_2}2 A^{\alpha -1}.
\]
Thus, the conclusion follows.
\end{proof}

From Proposition~\ref{P:DeltaN} and Proposition~\ref{P:squareroot},
we obtain
\begin{theorem} \label{T:DeltaN}
There exist $C>0$, which depends only on the set $\Lambda$ such that
if $s>\frac{n+3}2$
\[
|(-\Delta_{\p \Omega})^{\frac12} - \CN|_{L(H^{s'} (\p \Omega))} \le
C, \qquad s' \in [1-s, s-1]
\]
and if $s>2$ and $s \in (\frac {n+1}2, \frac {n+3}2)$, for $\alpha=
\frac {n+5}2 -s$,
\begin{equation}\label{E:deltacn4}
|(-\Delta_{\p \Omega})^{\frac12} - \CN|_{L(H^{s'} (\p \Omega),
H^{s'-\alpha+1} (\p \Omega))} \le C, \qquad s' \in (1-s,
\frac{n+1}2).
\end{equation}
\end{theorem}

\begin{proof}
We will give the proof for the second case only as the proof for the
first proof is similar. The estimate~\eqref{E:deltacn4} follows
directly from Proposition~\ref{P:DeltaN} and~\ref{P:squareroot} for
$s' \in (1-s, \frac{n-1}2)$. To prove the estimate for $s' \in
[\frac{n-1}2, \frac{n+1}2)$, we observe that $s'-2 \in (1-s,
\frac{n-1}2)$ and we have
\[
|(I - \Delta_{\p \Omega})^{-1} ( (-\Delta_{\p \Omega})^{\frac 12} -
\CN)(I - \Delta_{\p \Omega})|_{L(H^{\frac{n+1}2} (\p \Omega),
H^{\frac{n+1}2 -\alpha +1} (\p \Omega))} \le C.
\]
Thus~\eqref{E:deltacn4} follows from the commutator
estimate~\eqref{E:deltacn3}.
\end{proof}

\noindent {\bf Decomposition of vector fields.} We conclude this
section by introducing the velocity field decomposition. Given an
$L^2$ vector field $u: \Omega \to \RR^n$, it is standard to
decompose it into the divergence free part $v \in L^2$ and the
gradient part $-\nabla p$ for $p \in H^1_0(\Omega)$. In fact,
\begin{equation} \label{E:decomp1}
-\Delta p = \nabla \cdot u \qquad \qquad v = u + \nabla p.
\end{equation}
For any divergence free vector field $v \in L^2 (\Omega)$, the
normal component on the boundary $v^\perp \triangleq v\cdot N : \p
\Omega \to \RR$ in $H^{-\frac 12} (\p \Omega)$ is defined as
\[
< v^\perp, \psi> = \int_\Omega v\cdot \nabla \psi_\CH dx
\]
for any $\psi \in H^{\frac 12} (\p \Omega)$. By interpolation, for
any $s_1 \in [0, s-\frac 12]$,
\begin{equation} \label{E:vperp}
|v^\perp|_{H^{s_1-\frac 12} (\p \Omega)} \le C |v|_{H^{s_1}
(\Omega)},
\end{equation}
with $C$ uniform in $\Omega \in \Lambda$. This induces a
decomposition of $v$ into two divergence free parts, the rotation
part $v_r$ and the irrotational (or gradient) part $v_{ir}$, as
follows
\begin{equation} \label{E:decom2}
v_{ir} = \nabla \CH \CN^{-1} v^\perp, \qquad \qquad v_r = v- v_{ir}.
\end{equation}
It is easy to verify $v_r, \; v_{ir} \in L^2 (\Omega)$ and
\[
\nabla \cdot v_r = \nabla \cdot v_{ir}=0, \quad <v_r, v_{ir}>=0,
\quad v_r^\perp=0.
\]
If $v$ is a divergence free velocity field, $v_r$ component is
responsible of the internal rotation and $v_{ir}$ of the motion of
the domain.

\bigskip

\noindent{\bf  Notation}

\noindent
 \\  $\text{tr}(A)$: the trace of an operator.
 \\  $A^*$: the adjoint operator of an operator.
 \\  $A_1 \cdot A_2 =\text{tr} (A_1 (A_2)^*)$, for two operators.
 \\  $B(S, \ep) = \cup_{x \in S} B(x, \ep)$: an $\ep$-neighborhood of a set $S$.
 \\  $D$ and $\p$: differentiation with respect to spatial variables.
 \\  $\nabla f$: the gradient vector of a scalar function $f$.
 \\  $\nabla_X$:  the directional directive in the direction $X$.
 \\  $\perp$ and $\top$: the normal and the tangential components of the relevant quantities.
 \\  $\D_t = \p_t + v^i\p_{x^i}$ the material derivative  along the particle path.
 \\  $\D_t^\top$: the projection of $\D_t$ to the tangent space of $\p \Omega_t \subset \RR^n$.
\\  $N(t,x)$: the outward unit normal vector of $\p \Omega_t$ at $x \in \p \Omega_t$.
 \\  $\Pi$: the second fundamental form of $\p \Omega_t$, $\Pi(t, x)(w) = \nabla_w N \in T_x \p\Omega_t$.
\\ $\Pi(X, Y)=\Pi(X) \cdot Y$.
\\  $\kappa$: the mean curvature of $\p \Omega_t$, i.e. $\kappa = \text{tr} \Pi$.
\\$f_\CH=\CH(f)$: the harmonic extension of $f$ on $\Omega_t$.
\\$\CN (f) = \nabla_N \CH(f) : \p \Omega \to \RR$:  the Dirichlet-Neumann operator.
\\$\bar X = X\circ u^{-1}$ the Lagrangian coordinates description of $X$.
\\  $\CD$: the covariant differentiation on $ \p \Omega_t \subset \RR^n$.
\\$\CD_w = \nabla_w^\top$,  for any $ x \in \p\Omega_t \quad w\in T_x \p\Omega_t$.
 \\  $\R(X, Y)$, $\; X, Y\in T_x  \p\Omega_t$: the curvature tensor of $ \p\Omega_t$.
 \\  $\Delta_\CM \triangleq \text{tr} \CD^2$: the Beltrami-Lapalace operator on a Riemannian manifold $\CM$.
 \\$\Delta^{-1}$: the inverse Laplacian with zero Dirichlet data.
 \\$\Gamma = \{\phi \, : \Omega_t \to \RR^n \text{ ; volume preserving homeomorphism} \}$
 \\  $\bar \SD$: the covariant derivative on $\Gamma$,
 \\  $ \SD$: represent $\bar\SD$ in Eulerian coordinates.
 \\  $\bar \SR$: the curvature operator on $\Gamma$.
 \\  $ \SR$: represent $\bar\SR$ in Eulerian coordinates.
 \\  $\text{II}$: the second fundamental form of $\Gamma \subset L^2$
\\$\text{II}_u(w_1, w_2) = \nabla^\perp_{w_1}w_2$,  for any $u \in \Gamma, \quad w_1, w_2 \in T_u \Gamma$
\\ $p_{v,w} =-\Delta^{-1} \text{tr} (DvDw)$.

\end{document}